# The new Fermat-type factorization algorithm based on Fermat's factorization algorithm and the author's modified Euler's totient function

Savino Detto

Via Paolo Borsellino, 21/2 - 76012 - Canosa di Puglia (BT), Italy
*Email address*: **savinodetto@libero.it**

## Abstract

Let $n$ be any odd natural number other than a perfect square. We show that the new factorization algorithm, presented in this paper and which we call *DFM-1* (where DFM stands for Detto's Factorization Method), is much more efficient than the implementation technique of Fermat's Factorization Algorithm (FFA) called FFA-1, which, among the implementation techniques of Fermat's Factorization Algorithm (FFA), is the one that requires the fewest iterations to identify the non-trivial and trivial factors of $n$ (excluding the cases in which the two factors of the pair of non-trivial or trivial factors of $n$ are so close to each other that they can be identified at the $1^{st}$ iteration with each of such implementation techniques). Indeed, by the way in which Euler's totient function of any $n$ that is a semiprime is applied to FFA-1, we arrive at the new factorization algorithm (DFM-1), which halves (possibly rounding up to the next integer) the number of iterations required by FFA-1. Furthermore, in this paper, we present the *hypothetical scenario* according to which the number of iterations could possibly be further reduced. Finally, and still in relation to this new factorization algorithm, in this paper we present the *limit number of iterations*, which is less than the number of iterations required by DFM-1 to reach the condition $x - y = 1$ which characterizes the pair of trivial factors of $n$, beyond which it is no longer possible for pairs of non-trivial factors of $n$ to occur.

## 1. Introduction

The Fundamental Theorem of Arithmetic states that every natural number greater than 1 is either a prime or can be expressed as a product of primes. In the second case, the representation is unique, regardless of the order in which the prime factors appear.

In Number Theory, factorization means the decomposition of a composite integer into a set of non-trivial divisors (i.e. other than 1 and the composite integer itself that is being factored) which, multiplied together, determine the starting number. The non-trivial divisors are not necessarily prime (prime factors), but can also be, in whole or in part, composite (composite factors).

Nowadays, several factorization methods or algorithms are known, some are simple and others are more complex and refined, both deterministic and probabilistic such as (for example, to indicate just a few) the algorithms described in [5, pp.132-151], and, in general, their efficiency is always based on the characteristics of the prime factors; but, currently, a polynomial solution for the factorization is not yet known. And, probably, a polynomial solution for the factorization may not exist. Indeed, many modern cybersecurity applications (cryptographic applications) are based precisely on the assumption that a polynomial solution for factorization does not exist. To tell the truth, to date it has not been proven that this solution exists, but neither it has been proven that it cannot exist; so, currently the search for this polynomial solution for factorization remains a still open problem in Number Theory.

Thus, Fermat's Factorization Algorithm (FFA) [1, p. 6; 3, p. 1470], discovered by the French mathematician and magistrate Pierre De Fermat in 1600, is one of the factorization algorithms known today. It is more efficient the smaller the difference between two non-trivial factors, prime or composite, of $n$, with $n$ being any odd natural number other than a perfect square. This factorization algorithm immediately captures our attention and arouses our interest, leading us to investigate it in a particularly detailed and thorough manner. We hope to be able, one day, to publish the in-depth studies on Fermat's Factorization Algorithm (FFA); thus, the reader will learn, inter alia, about the path that, in 2019, led us to distinguish all odd natural numbers (therefore, including any odd perfect square and, therefore, including the number 1 which is a perfect square) into two different categories that are easily and immediately identifiable and one of which excludes, with absolute certainty, the odd natural numbers which are perfect squares.

Following in-depth studies into several distinct topics covered in a study on natural numbers, and focusing attention on the implementation technique of Fermat's Factorization Algorithm (FFA) called FFA-1 [2; 3, p. 1470], of which we remain aware that it is not the most efficient factorization algorithm currently known (if we exclude the circumstance in which the two factors of the pair of non-trivial or trivial factors of $n$, with $n$ being any odd natural number other than a perfect square, are so close to each other that they can be identified at the 1st iteration with FFA-1, and this happens when the sum of the two factors of the pair of non-trivial or trivial factors of $n$ is equal to $2 \cdot \lceil\sqrt{n}\rceil$, i.e. when the difference between the two factors of the pair of non-trivial or trivial factors of $n$ is less than or equal to $\lfloor\sqrt{2^3} \cdot \sqrt[4]{n}\rfloor$), we arrived at the new factorization algorithm, which is presented in this paper and which we call *DFM-1* (where DFM stands for Detto's Factorization Method), whose efficiency will prove to be clearly and reasonably greater than that of FFA-1, which, among the implementation techniques of Fermat's Factorization Algorithm (FFA), is the one that requires the fewest iterations to identify the non-trivial and trivial factors of $n$ (excluding the cases in which the two factors of the pair of non-trivial or trivial factors of $n$ are so close to each other that they can be identified at the 1st iteration with each of such implementation techniques). Briefly, we start by considering the unequivocal and absolutely evident fact that the value of Euler's totient function $\varphi(n)$, which is associated with any $n$ that is a semiprime (with $n$ being any odd natural number other than a perfect square), is a positive integer that is certainly always a multiple of 4 and, therefore, it is certainly always divisible by 4 (obviously, also the value of Euler's totient function $\varphi(n)$, which is associated with any $n$ that is a non-semiprime composite formed of at least two distinct primes, is a positive integer that is certainly always a multiple of 4 and, therefore, it is certainly always divisible by 4); this gradually leads us to make other considerations and observations, putting a particular light on Euler's totient function (by modifying, for factorization purposes, Euler's totient function $\varphi(n)$ of any $n$ that is a semiprime), with the related in-depth studies that lead us first to the new factorization algorithm (DFM-1), of which we nonetheless always remain aware that it is not the most efficient factorization algorithm (if we exclude the cases, which are infinite, in which the two factors of the pair of non-trivial or trivial factors of $n$ are so close to each other that they can be identified at the 1st iteration with DFM-1, and this happens when the sum of the two factors of the pair of non-trivial or trivial factors of $n$ is equal to $n - \left\lfloor \frac{n - 2 \cdot \lfloor\sqrt{n}\rfloor}{4} \right\rfloor \cdot 4 + 1$, i.e. when the difference between the two factors of the pair of non-trivial or trivial factors of $n$ is less than or equal to the nearest integer to $2^2 \cdot \sqrt[4]{n}$), and then to the *hypothetical scenario*, with which further and interesting perspectives on factorization open up, and to the representation of the *limit number of iterations*—which is less than the number of iterations required by DFM-1 to reach the condition $x - y = 1$ which characterizes the pair of trivial factors of $n$—beyond which it is no longer possible for pairs of non-trivial factors of $n$ to occur.

In summary, the new factorization algorithm (DFM-1), presented in this paper, puts an end to the step size of 1, which is adopted with FFA-1 [2; 3, p.1470] and which is certainly replaced by the step size of 2, for increasing the value of $x$ in the equation $y = \sqrt{x^2 - n}$ in order to identify the non-trivial factors of $n$ or to determine the primality of $n$, with $n$ being any odd natural number other than a perfect square. Having modified, for factorization purposes, Euler's totient function $\varphi(n)$ of any $n$ that is a semiprime (i.e. having modified, for factorization purposes, Euler's totient function $\varphi(n)$ of any odd semiprime other than a perfect square), it will suffice to calculate the initial value of $x$ as described in this paper to proceed with a step size of 2, rather than 1, when incrementing the value of $x$ in the equation $y = \sqrt{x^2 - n}$ in order to identify the non-trivial factors of $n$ or to determine the primality of $n$.

Therefore, the reasons why, in this paper, we *will compare* the *new factorization algorithm* (DFM-1), which is presented in this paper, *solely* with FFA-1 [2; 3, p. 1470] and not with other *factorization algorithms* are already clear and reasonable. However, remaining always aware that, currently, there are *much more efficient factorization algorithms* than DFM-1; nevertheless, uncovering new results in mathematics and Number Theory is always valuable, and knowledge benefits from it.

## 2. Fermat's Factorization Algorithm (FFA)

Fermat's Factorization Algorithm (FFA) [1, p. 6; 3, p. 1470] is based on the fact that every odd natural number, other than a perfect square (and, therefore, also other than 1 which is a perfect square), can be expressed as a difference of two perfect squares, as represented by the equation indicated below, with $n$ being any odd natural number other than a perfect square and with $x$ and $y$ which are positive integers:

$$n = x^2 - y^2, \qquad (2.1a)$$

i.e.:

$$n = (x + y) \cdot (x - y). \tag{2.1b}$$

Obviously, any odd natural number that is a perfect square (but other than 1, since $x$ and $y$ are positive integers and, therefore, with $y \neq 0$) can also be represented by the two equations indicated above, but only with reference to the pair of trivial factors (as in the case of primes) if only one pair of non-trivial factors occurs, and, if more than one pair of non-trivial factors occurs, also with reference to any pair of non-trivial factors that is subsequent to the first pair of non-trivial factors, since, obviously, with the first or unique pair of non-trivial factors of a perfect square —the two factors being equal to each other— we have $y = 0$, with $y$, therefore, which is not a positive integer; but, by convention (and convenience), in this paper we exclude a priori that $n$ is a perfect square.

That said, and denoting by $a$ and $b$, with $a > b$, the factors, prime or composite and also with $a = n$ and $b = 1$, of any pair of non-trivial factors of $n$ and of the pair of trivial factors of $n$ (obviously, $a$ and $b$ are odd), we have:

$$n = a \cdot b, \tag{2.2a}$$

with this last equation which can also be written as follows:

$$n = \left(\frac{a + b}{2}\right)^2 - \left(\frac{a - b}{2}\right)^2 \tag{2.2b}$$

whence, in view of (2.1$a$), we deduce that:

$$x = \frac{a + b}{2} \tag{2.3}$$

and

$$y = \frac{a - b}{2}, \tag{2.4}$$

and therefore, combining (2.1$a$) and (2.2$a$), we have:

$$n = x^2 - y^2 = a \cdot b \tag{2.5}$$

since we have:

$$
\begin{aligned}
n &= x^2 - y^2 \\
&= \left(\frac{a+b}{2}\right)^2 - \left(\frac{a-b}{2}\right)^2 \\
&= \frac{(a+b)^2}{4} - \frac{(a-b)^2}{4} \\
&= \frac{a^2 + 2 \cdot a \cdot b + b^2 - (a^2 - 2 \cdot a \cdot b + b^2)}{4} \\
&= \frac{a^2 + 2 \cdot a \cdot b + b^2 - a^2 + 2 \cdot a \cdot b - b^2}{4} \\
&= \frac{4 \cdot a \cdot b}{4} \\
&= a \cdot b,
\end{aligned}
$$

so arriving to represent the following equation:

$$x^2 - y^2 - n = 0 \tag{2.6}$$

and, consequently, to represent the following equation:

$$n = x^2 - y^2, \tag{2.1$a$}$$

which is emblematic of Fermat's Factorization Algorithm (FFA) and which can also be written as follows:

$$y = \sqrt{x^2 - n}. \tag{2.7}$$

In practice, searching for the pair of non-trivial factors $a$ and $b$ of $n$ means searching for the pair of positive integers $x$ and $y$ (with $x - y \neq 1$, since $x - y = 1$ is the condition that characterizes the pair of trivial factors of $n$) such that the sum $S$ vanishes:

$$S = x^2 - y^2 - n = 0. \tag{2.8}$$

Therefore, having found the pair of positive integers $x$ and $y$ with $x - y \neq 1$, we can proceed to determine the non-trivial factors $a$ and $b$ of $n$ as shown below:

$$a = x + y \tag{2.9}$$

and

$$b = x - y \tag{2.10}$$

since, from (2.3) and (2.4), we have , respectively:

$$a = 2 \cdot x - b \tag{2.11}$$

and

$$b = a - 2 \cdot y \tag{2.12}$$

and, therefore, we have , respectively:

$$\begin{aligned} a &= 2 \cdot x - b \\ &= 2 \cdot x - (a - 2 \cdot y) \\ &= 2 \cdot x - a + 2 \cdot y \end{aligned}$$

from which it follows that:

$$2 \cdot a = 2 \cdot (x + y)$$

which simplifies to:

$$a = x + y, \tag{2.9}$$

and

$$\begin{aligned} b &= a - 2 \cdot y \\ &= (2 \cdot x - b) - 2 \cdot y \\ &= 2 \cdot x - b - 2 \cdot y \end{aligned}$$

from which it follows that:

$$2 \cdot b = 2 \cdot (x - y)$$

which simplifies to:

$$b = x - y. \tag{2.10}$$

## 3. Euler's Totient Function

In short, Euler's totient function $\varphi(N)$ [4] is a function which represents, for any positive integer $N$, the number of positive integers less than $N$ and that are coprime with $N$.

The equation that defines the value of Euler's totient function $\varphi(N)$ can be written as follows:

$$\varphi(N) = \varphi(p_1^{k_1}) \cdot \varphi(p_2^{k_2}) \cdot \ \ldots \ \cdot \varphi(p_f^{k_f}), \tag{3.1$a$}$$

i.e.:

$$\varphi(N) = \left[(p_1 - 1) \cdot p_1^{k_1 - 1}\right] \cdot \left[(p_2 - 1) \cdot p_2^{k_2 - 1}\right] \cdot \ \dots \ \cdot \left[(p_f - 1) \cdot p_f^{k_f - 1}\right], \quad (3.1b)$$

where $p_1, p_2, \dots, p_f$ are all the distinct prime factors, i.e. prime factors that are different from each other, which form $N$, and $k_1, k_2, \dots, k_f$ are the exponents of the respective distinct prime factors that form $N$.

Euler's totient function $\varphi(N)$ is multiplicative between two positive integers $a$ and $b$, but only if $a$ and $b$ are prime to each other, i.e. they are coprime such that $\text{GCD}(a, b) = 1$. Indeed, if we consider, as an example, the positive integers $a = 35 = 7 \cdot 5$ and $b = 33 = 11 \cdot 3$, with $a$ and $b$ as factors of the pair of non-trivial factors of $N = 1155$, we have:

$$\begin{aligned} \varphi(N) = \varphi(1155) &= \varphi(a) \cdot \varphi(b) \\ &= \varphi(35) \cdot \varphi(33) \\ &= 24 \cdot 20 \\ &= 480 \end{aligned}$$

since, with $p_2 = 5$, $p_3 = 7$, $k_2 = 1$, and $k_3 = 1$, we have:

$$\begin{aligned} \varphi(a) &= \varphi(p_2^{k_2}) \cdot \varphi(p_3^{k_3}) \\ &= \left[(p_2 - 1) \cdot p_2^{k_2 - 1}\right] \cdot \left[(p_3 - 1) \cdot p_3^{k_3 - 1}\right], \end{aligned}$$

and substituting their respective values, we obtain:

$$\begin{aligned} \varphi(35) &= \varphi(5^1) \cdot \varphi(7^1) \\ &= \left[(5 - 1) \cdot 5^{1 - 1}\right] \cdot \left[(7 - 1) \cdot 7^{1 - 1}\right] \\ &= 4 \cdot 6 \\ &= 24, \end{aligned}$$

and since, with $p_1 = 3$, $p_4 = 11$, $k_1 = 1$, and $k_4 = 1$, we have:

$$\begin{aligned} \varphi(b) &= \varphi(p_1^{k_1}) \cdot \varphi(p_4^{k_4}) \\ &= \left[(p_1 - 1) \cdot p_1^{k_1 - 1}\right] \cdot \left[(p_4 - 1) \cdot p_4^{k_4 - 1}\right], \end{aligned}$$

and substituting their respective values, we obtain:

$$\begin{aligned} \varphi(33) &= \varphi(3^1) \cdot \varphi(11^1) \\ &= \left[(3 - 1) \cdot 3^{1 - 1}\right] \cdot \left[(11 - 1) \cdot 11^{1 - 1}\right] \\ &= 2 \cdot 10 \\ &= 20, \end{aligned}$$

with the value of $\varphi(N) = \varphi(1155) = 480$, which is $\varphi(a) \cdot \varphi(b) = \varphi(35) \cdot \varphi(33) = 24 \cdot 20$, which coincides, since $a$ and $b$ are coprime (since $a$ and $b$ have no prime factors in common), with the correct value of $\varphi(N) = \varphi(1155) = 480$ as calculated below, with $p_1 = 3$, $p_2 = 5, p_3 = 7, p_4 = 11, k_1 = 1, k_2 = 1, k_3 = 1$, and $k_4 = 1$:

$$\begin{aligned}\varphi(N) &= \varphi(p_1^{k_1}) \cdot \varphi(p_2^{k_2}) \cdot \varphi(p_3^{k_3}) \cdot \varphi(p_4^{k_4}) \\ &= \left[(p_1 - 1) \cdot p_1^{k_1 - 1}\right] \cdot \left[(p_2 - 1) \cdot p_2^{k_2 - 1}\right] \cdot \left[(p_3 - 1) \cdot p_3^{k_3 - 1}\right] \cdot \left[(p_4 - 1) \cdot p_4^{k_4 - 1}\right],\end{aligned}$$

and substituting their respective values, we obtain:

$$\begin{aligned}\varphi(1155) &= \varphi(3^1) \cdot \varphi(5^1) \cdot \varphi(7^1) \cdot \varphi(11^1) \\ &= \left[(3 - 1) \cdot 3^{1 - 1}\right] \cdot \left[(5 - 1) \cdot 5^{1 - 1}\right] \cdot \left[(7 - 1) \cdot 7^{1 - 1}\right] \cdot \left[(11 - 1) \cdot 11^{1 - 1}\right] \\ &= 2 \cdot 4 \cdot 6 \cdot 10 \\ &= 480.\end{aligned}$$

While, if we consider, as another example, the positive integers $a = 30 = 2 \cdot 3 \cdot 5$ and $b = 21 = 3 \cdot 7$, with $a$ and $b$ as factors of the pair of non-trivial factors of $N = 630$, we have:

$$\begin{aligned}\varphi(N) = \varphi(630) &= \varphi(a) \cdot \varphi(b) \\ &= \varphi(30) \cdot \varphi(21) \\ &= 8 \cdot 12 \\ &= 96\end{aligned}$$

since, with $p_1 = 2, p_2 = 3, p_3 = 5, k_1 = 1, k_2 = 1$, and $k_3 = 1$, we have:

$$\begin{aligned}\varphi(a) &= \varphi(p_1^{k_1}) \cdot \varphi(p_2^{k_2}) \cdot \varphi(p_3^{k_3}) \\ &= \left[(p_1 - 1) \cdot p_1^{k_1 - 1}\right] \cdot \left[(p_2 - 1) \cdot p_2^{k_2 - 1}\right] \cdot \left[(p_3 - 1) \cdot p_3^{k_3 - 1}\right],\end{aligned}$$

and substituting their respective values, we obtain:

$$\begin{aligned}\varphi(30) &= \varphi(2^1) \cdot \varphi(3^1) \cdot \varphi(5^1) \\ &= \left[(2 - 1) \cdot 2^{1 - 1}\right] \cdot \left[(3 - 1) \cdot 3^{1 - 1}\right] \cdot \left[(5 - 1) \cdot 5^{1 - 1}\right] \\ &= 1 \cdot 2 \cdot 4 \\ &= 8,\end{aligned}$$

and since, with $p_2 = 3, p_4 = 7, k_2 = 1$, and $k_4 = 1$, we have:

$$\begin{aligned}\varphi(b) &= \varphi(p_2^{k_2}) \cdot \varphi(p_4^{k_4}) \\ &= \left[(p_2 - 1) \cdot p_2^{k_2 - 1}\right] \cdot \left[(p_4 - 1) \cdot p_4^{k_4 - 1}\right],\end{aligned}$$

and substituting their respective values, we obtain:

$$\begin{aligned}\varphi(21) &= \varphi(3^1) \cdot \varphi(7^1)\\ &= \left[(3-1) \cdot 3^{1-1}\right] \cdot \left[(7-1) \cdot 7^{1-1}\right]\\ &= 2 \cdot 6\\ &= 12,\end{aligned}$$

with the value of $\varphi(N) = \varphi(630) = 96$, which is $\varphi(a) \cdot \varphi(b) = \varphi(30) \cdot \varphi(21) = 8 \cdot 12$, which does not coincide, since $a$ and $b$ are not coprime (since $a$ and $b$ have the prime 3 in common), with the correct value of $\varphi(N) = \varphi(630) = 144$ as calculated below, with $p_1 = 2$, $p_2 = 3, p_3 = 5, p_4 = 7, k_1 = 1, k_2 = 2, k_3 = 1$, and $k_4 = 1$:

$$\begin{aligned}\varphi(N) &= \varphi(p_1^{k_1}) \cdot \varphi(p_2^{k_2}) \cdot \varphi(p_3^{k_3}) \cdot \varphi(p_4^{k_4})\\ &= \left[(p_1-1) \cdot p_1^{k_1-1}\right] \cdot \left[(p_2-1) \cdot p_2^{k_2-1}\right] \cdot \left[(p_3-1) \cdot p_3^{k_3-1}\right] \cdot \left[(p_4-1) \cdot p_4^{k_4-1}\right],\end{aligned}$$

and substituting their respective values, we obtain:

$$\begin{aligned}\varphi(630) &= \varphi(2^1) \cdot \varphi(3^2) \cdot \varphi(5^1) \cdot \varphi(7^1)\\ &= \left[(2-1) \cdot 2^{1-1}\right] \cdot \left[(3-1) \cdot 3^{2-1}\right] \cdot \left[(5-1) \cdot 5^{1-1}\right] \cdot \left[(7-1) \cdot 7^{1-1}\right]\\ &= 1 \cdot 6 \cdot 4 \cdot 6\\ &= 144.\end{aligned}$$

That said, it must be emphasized right away that, for the purposes of applying the new factorization algorithm (DFM-1) presented in this paper, both factors $a$ and $b$, with $a > b$, of any possible pair of non-trivial factors of $n$ and of the pair of trivial factors of $n$ (with $n$ being any odd natural number other than a perfect square) will in any case always be considered as prime even if they are not (therefore, also both factors $a$ and $b$ of the pair of trivial factors of $n$, with $a$ which might not even be prime, will in any case always be considered as prime, just as also both factors $a$ and $b$ of any possible pair of non-trivial factors of $n$, which have prime factors in common and, therefore, are not coprime, will in any case always be considered as prime), thus putting a particular light on Euler's totient function (having modified, for factorization purposes, Euler's totient function $\varphi(n)$ of any $n$ that is a semiprime,), inserting a new value, which is $\varphi_s(n)$, which will replace the value of $\varphi(n)$ for any $n$ that is a semiprime and will coincide with the value of such $\varphi(n)$ only in the case in which the value of $\varphi_s(n)$ is associated with the unique pair of non-trivial factors of any $n$ that is a semiprime, i.e. of any odd semiprime other than a perfect square.

Therefore, for the purposes of applying the new factorization algorithm (DFM-1) presented in this paper, the following equation that defines the value of Euler's totient function $\varphi(n)$ for any $n$ that is a semiprime:

$$\varphi(n) = \left[(p_1-1) \cdot p_1^{k_1-1}\right] \cdot \left[(p_2-1) \cdot p_2^{k_2-1}\right], \qquad (3.2a)$$

i.e.:

$$\varphi(n) = (p_1 - 1) \cdot (p_2 - 1), \tag{3.2b}$$

will be replaced by the following equation that defines the value of the product $(a - 1) \cdot (b - 1)$, which we denote by $\varphi_s(n)$:

$$\varphi_s(n) = (a - 1) \cdot (b - 1), \tag{3.3}$$

which is associated with any possible pair of non-trivial factors *a* and *b* of *n* and to the pair of trivial factors *a* and *b* of *n*, with *a* and *b*, and with $a > b$, which are, therefore, the factors of any possible pair of non-trivial factors of *n* and of the pair of trivial factors of *n*, which will in any case always be considered as prime even if they are not, and which may have prime factors in common and, thus, may not be coprime. Thus, these are all the aspects representing the modification, for factoring purposes, of Euler's totient function $\varphi(n)$ of any *n* that is a semiprime, i.e. of any odd semiprime other than a perfect square.

## 4. FFA-1 Implementation Technique

As already highlighted in the Abstract and the Introduction of this paper, FFA-1 [2; 3, p. 1470], which is referred to in this paper for the purposes of comparison with the new factorization algorithm (DFM-1), is the implementation technique, among the implementation techniques [1, pp. 6-8; 2; 3, pp. 1470-1471] of Fermat's Factorization Algorithm (FFA), which requires the fewest iterations to identify the non-trivial and trivial factors of *n* (excluding the cases in which the two factors of the pair of non-trivial or trivial factors of *n* are so close to each other that they can be identified at the 1st iteration with each of such implementation techniques), with *n* being any odd natural number other than a perfect square.

Before starting to describe FFA-1, it is appropriate to highlight that there are two other implementation techniques of Fermat's Factorization Algorithm (FFA), which differ from FFA-1 also due to the absence of the square root operation and for which we will not proceed with their respective descriptions as they are not topics of discussion in this paper. These two additional implementation techniques are the implementation technique—which we call FFA-2*a*—described in [1, pp. 6-8], and the implementation technique—which we call FFA-2*b* here and which is already called FFA-2 in [3, p. 1471]—described in [3, p. 1471]. Practically, FFA-2*a* and FFA-2*b* (i.e. FFA-2) are the same implementation technique, but it is important to highlight that FFA-2*a* is more efficient than FFA-2*b* (i.e. FFA-2) since with the latter, in addition to having to perform more operations overall compared to FFA-2*a*, larger numbers are handled, for the same iteration, compared to those handled with FFA-2*a*, and, therefore and for the same *n* considered, FFA-2*b* (i.e. FFA-2) has a higher computational cost than FFA-2*a*.

That said, we will now describe FFA-1.

The FFA-1 implementation technique consists in finding, by repeatedly applying the following equation already given in Section 2:

$$y = \sqrt{x^2 - n}, \tag{2.7}$$

the possible unique pair of positive integers $x$ and $y$ (with $x - y \neq 1$) or the possible different pairs of positive integers $x$ and $y$ (with $x - y \neq 1$) such as to be able to identify the possible unique pair of non-trivial factors $a$ and $b$ of $n$ or the possible different pairs of non-trivial factors $a$ and $b$ of $n$, with $a > b$, starting from the following value of $x = x_1$:

$$x_1 = \lceil \sqrt{n} \rceil, \tag{4.1}$$

which is the smallest integer greater than $\sqrt{n}$, and proceeding, thereafter, always with a step size of 1 for the increase in the value of $x$, i.e. increasing the value of $x$ always by one unit at a time (starting, therefore, from the 2nd iteration with the value of $x_2 = x_1 + 1 = \lceil \sqrt{n} \rceil + 1$), until we are certain that pairs of positive integers $x$ and $y$ (with $x - y \neq 1$) can no longer occur such that we can identify pairs of non-trivial factors $a$ and $b$ of $n$ in addition to the pairs of non-trivial factors already identified, or, possibly, until we are certain of the primality of $n$. It is necessary to keep in mind that, by construction, we have that $x \geq \sqrt{n}$, with $x$ being an integer; therefore, it will be necessary to check, preliminarily, that $n$ is not a perfect square, i.e. that the condition $x = \sqrt{n}$ does not occur. Furthermore, before proceeding with FFA-1, we will preliminarily subject $n$ to a quick probabilistic test of primality such as [5, pp. 258-259]. Obviously, with FFA-1 we will always have, for each iteration, $S = 0$—see (2.8)—which therefore holds even when the value of $y$ is not an integer.

It should be highlighted that, in the cases of $n = 5$ and $n = 3$, the value of $x_1 = \lceil \sqrt{n} \rceil$ coincides with the value of $x = \frac{n+1}{2}$ with which the pair of trivial factors of n is identified.

Therefore, for each pair of positive integers $x$ and $y$ (with $x - y \neq 1$) that may be found, we will proceed to determine the non-trivial factors $a$ and $b$ of $n$ with the following equations already given in Section 2, highlighting that, with $n$ for which more than one pair of non-trivial factors occurs, the first identified pair of non-trivial factors of $n$ is the one with the smallest non-trivial value of $a$ and the largest non-trivial value of $b$, i.e. with non-trivial $a$ and $b$ closest to each other; while, the last identified pair of non-trivial factors of $n$ is the one with the largest non-trivial value of $a$ and the smallest non-trivial value of $b$, i.e. with non-trivial $a$ and $b$ farthest from each other:

$$a = x + y \tag{2.9}$$

and

$$b = x - y. \tag{2.10}$$

Below, we proceed with a simple example of the application of FFA-1, however highlighting (as well pointed out in the Introduction of the book, in Italian, *I moltiplicatori complessi per la possibile scomposizione dei Numeri RSA* [6, p. 22]) that it is advisable, since it is more convenient, to stop at the identification of the first (which could also be the only one) pair of positive integers $x$ and $y$ with $x - y \neq 1$ (but as long as there is a pair of positive integers $x$ and $y$ with $x - y \neq 1$), proceeding, obviously with a step size of 1 for the increase in the value of $x$, to check the factorization (or possible primality), by adopting FFA-1 again, of the first (or unique, but we cannot know a priori) identified non-trivial factors $a$ and $b$ of $n$ (obviously, after checking that $a$ and $b$ are not perfect squares and after having preliminarily subjected $a$ and $b$ to a quick probabilistic test of primality). Furthermore, it is also advisable to resort to the Trial Division Algorithm (TDA) [5, pp. 133-134], before proceeding with FFA-1, attempting to divide $n$ by a relatively small number (relative to the order of magnitude of $n$) of successive primes, starting from the number 3, contained in $\sqrt{n}$ and possibly finding prime factors of small magnitude so as to possibly have to proceed to check the factorization (or possible primality) of a $n$ smaller than the one initially considered.

Let $n = 70399$, which is the product of the primes 113, 89 and 7, be the composite to factorize. After checking that $n$ is not a perfect square, i.e. that the condition $x = \sqrt{70399}$ does not occur, and after having preliminarily subjected $n$ to a quick probabilistic test of primality, we proceed, by repeatedly applying of (2.7) (i.e. the equation $y = \sqrt{x^2 - n}$), with the search for the pairs of non-trivial factors $a$ and $b$ of $n$ until reaching the condition $x - y = 1$, starting from the value of $x = x_1 = \lceil\sqrt{n}\rceil = 266$, which is the smallest integer greater than $\sqrt{n} = 265{,}32809\ldots$, and proceeding, thereafter, always with a step size of 1 for the increase in the value of $x$, with the results reported in the following Table 1 where the number of iterations is denoted by $i$:

Table 1. Example of the application of FFA-1 with $n = 70399$

| $i$ | $x$ | $y$ | $S$ | $a$ | $b$ |
|---|---|---|---|---|---|
| 1 | 266 | non-integer | 0 | ----- | ----- |
| 2 | 267 | non-integer | 0 | ----- | ----- |
| 3 | 268 | non-integer | 0 | ----- | ----- |
| 4 | 269 | non-integer | 0 | ----- | ----- |
| OMISSED | | | | | |
| 102 | 367 | non-integer | 0 | ----- | ----- |
| **103** | **368** | **255** | **0** | **623** | **113** |
| 104 | 369 | non-integer | 0 | ----- | ----- |
| 105 | 370 | non-integer | 0 | ----- | ----- |
| 106 | 371 | non-integer | 0 | ----- | ----- |
| OMISSED | | | | | |
| 174 | 439 | non-integer | 0 | ----- | ----- |
| **175** | **440** | **351** | **0** | **791** | **89** |
| 176 | 441 | non-integer | 0 | ----- | ----- |
| 177 | 442 | non-integer | 0 | ----- | ----- |
| 178 | 443 | non-integer | 0 | ----- | ----- |
| OMISSED | | | | | |
| 4766 | 5031 | non-integer | 0 | ----- | ----- |
| **4767** | **5032** | **5025** | **0** | **10057** | **7** |
| 4768 | 5033 | non-integer | 0 | ----- | ----- |
| 4769 | 5034 | non-integer | 0 | ----- | ----- |
| OMISSED | | | | | |
| 34934 | 35199 | non-integer | 0 | ----- | ----- |
| **34935** | **35200** | **35199** | **0** | **70399** | **1** |

As we can easily deduce from Table 1 and denoting by $x_i$ the value of $x$ corresponding to a specific iteration considered, with FFA-1 we have:

$$i = x_i - \lfloor\sqrt{n}\rfloor \tag{4.2}$$

for each iteration until the identification of the pair of trivial factors of $n$; therefore, with FFA-1 we always have:

$$i = \frac{a+b}{2} - \lfloor\sqrt{n}\rfloor \tag{4.3}$$

both for the identification of any pair of non-trivial factors of $n$ and for the identification of the pair of trivial factors of $n$, and, in the latter case, we can also write the following equation:

$$i = \frac{n+1}{2} - \lfloor\sqrt{n}\rfloor, \tag{4.4}$$

since $a = n$ and $b = 1$, with the value of $i$ which, in the latter case, we can obviously know a priori.

That said, it is appropriate to highlight the circumstance (which is already known, as indicated in reference [2], and which is anyway easy to see) that, after the 1st iteration, the subsequent iterations can also be performed differently, i.e. without calculating the square of $x$ each time and without calculating the difference between $x^2$ and $n$ each time, thus reducing the computational cost (but always leaving the number of iterations to be performed unchanged). Indeed, starting from considering the known value (obviously integer) of $y^2$ corresponding to the 1st iteration, and denoting by $x_f$ the value of $x$ corresponding to the previous iteration and by $y_f^2$ and $y_r^2$, respectively, the value of $y^2 = x^2 - n$ corresponding to the previous iteration and the value of $y^2 = x^2 - n$ corresponding to the subsequent iteration, with FFA-1 we can also proceed by adopting the following equation:

$$y_r^2 = y_f^2 + 2 \cdot x_f + 1, \tag{4.5}$$

for each iteration subsequent to the 1st iteration and until the identification of the pair of trivial factors of $n$ (starting, therefore, from the 2nd iteration with the value of $x_f = x_1 = \lceil\sqrt{n}\rceil$ and, thereafter, increasing the value of $x_f$ always by one unit at a time until the identification of the pair of trivial factors of $n$), and, then, checking whether the value of the square root of $y_r^2$, corresponding to the subsequent iteration, is an integer.

Furthermore, it is possible to proceed with a further simplification (which is also already known, as indicated in reference [2]), reducing the number of square root operations to be performed and, therefore, further reducing the computational cost (but always leaving the number of iterations to be performed unchanged). Indeed, since, as we know, a perfect square can only end with specific digits (0, 1, 4, 5, 6, 9 in base 10), if the last digit of $y_1^2$ (where $y_1^2$ is the first value of $y_f^2$ possibly to be considered) or of $y_r^2$ (like the last digit of any $y^2 = x^2 - n$) is 2, 3, 7 or 8, we can eliminate the corresponding square root operation.

We conclude this section by highlighting that it would be possible to proceed with a further simplification of FFA-1, further reducing the computational cost (but always leaving the number of iterations to be performed unchanged). Indeed, and as already highlighted in the Introduction to this paper, following in-depth studies on Fermat's Factorization Algorithm (FFA), we have come to distinguish all odd natural numbers (therefore, including any odd perfect square and, therefore, including the number 1 which is a perfect square) into two different categories that are easily and immediately identifiable and one of which excludes, with absolute certainty, the odd natural numbers which are perfect squares. Briefly, these two categories are distinguished by a value that I called the discriminant $\delta$, encompassing the odd natural numbers characterized by $\delta = 3$ and the odd natural numbers characterized by $\delta = 5$, with the value of $\delta$ being easily and immediately identifiable. The odd natural numbers characterized by $\delta = 3$ include primes and composites other than perfect squares (indistinguishable from each other); while, the odd natural numbers characterized by $\delta = 5$ include primes, perfect squares (including the number 1 which is a perfect square) and composites other than perfect squares (indistinguishable from each other). Therefore, if the last digit of $y_1{}^2$ or of $y_r{}^2$ (like the last digit of any $y^2 = x^2 - n$) is 1, 5 or 9, before performing the corresponding square root operation we will proceed to check the corresponding discriminant $\delta$: if $y_1{}^2$ or $y_r{}^2$ (like any $y^2 = x^2 - n$) is characterized by $\delta = 3$, we will not perform the corresponding square root operation since, in such case, $y_1{}^2$ or $y_r{}^2$ (like any $y^2 = x^2 - n$) can only be a prime or a composite other than a perfect square, and, therefore, it certainly cannot be a perfect square; while, if $y_1{}^2$ or $y_r{}^2$ (like any $y^2 = x^2 - n$) is characterized by $\delta = 5$, we will necessarily have to perform the corresponding square root operation since, in such case, $y_1{}^2$ or $y_r{}^2$ (like any $y^2 = x^2 - n$) can be a prime, a perfect square or a composite other than a perfect square, and, therefore, it could be a perfect square. But, in this section (and, therefore, in this paper), we will not proceed to present this further simplification of FFA-1, just as we will not proceed to present the way to identify the value of the discriminant $\delta$ of any odd natural number, since, with the inevitable easy observations on this further simplification and the consequent considerations with the related in-depth studies both on FFA-1 and on the aforementioned distinction of all odd natural numbers into two different categories, one would inevitably and easily arrive at a further simplification—which was nevertheless noticed after the new factorization algorithm (DFM-1) presented in this paper—that is so substantial, as it reduces the number of iterations to be performed (obviously, excluding the cases in which the two factors of the pair of non-trivial or trivial factors of $n$ are so close to each other that they can be identified at the 1st iteration with FFA-1), that it determines, in practice, a further new factorization algorithm (which is, therefore, not yet known), with the same results as DFM-1 presented in this paper (but with different premises which, unlike DFM-1 presented in this paper, exclude the treatment of particular and important aspects which, instead, are presented in this paper, such as the modification, for factorization purposes, of Euler's totient function $\varphi(n)$ of any $n$ that is a semiprime, and as the hypothetical scenario), and, therefore, this substantial simplification, if presented now, would take us away from the topics covered in this paper, topics which are the new factorization algorithm (DFM-1)—with all its aspects—presented in this paper, and the comparison with FFA-1 as it is known today [2; 3, p. 1470]. Therefore, we postpone the presentation of this substantial simplification—together with both the presentation of the aforementioned simplification due to the discriminant $\delta$ and the presentation of the way to identify the value of the discriminant $\delta$ of any odd natural number—to a forthcoming paper, but already highlightyng the fact (which will be evident in such forthcoming paper) that the above-mentioned distinction (an easy and immediate distinction) of all odd natural numbers into two different categories, although it always remains valid, is useless both for such substantial simplification—with or without the respective simplifications that will be presented in such forthcoming paper—and for the new factorization algorithm (DFM-1)—with or without the simplifications that will be presented in the last part of Section 5 of this paper—presented in this paper.

## 5. Summary Presentation of DFM-1

The new factorization algorithm (DFM-1) in question, which is applicable—more efficiently than FFA-1 [2; 3, p. 1470]—to any odd natural number other than a perfect square, is still based on the following note (2.1*a*) that is the emblematic equation of Fermat's Factorization Algorithm (FFA), with $n$ being any odd natural number other than a perfect square and with $x$ and $y$ which are positive integers:

$$n = x^2 - y^2, \tag{2.1a}$$

but we will look for, by repeatedly applying of (2.7) (i.e. the equation $y = \sqrt{x^2 - n}$), the possible unique pair of positive integers $x$ and $y$ (with $x - y \neq 1$) or the possible different pairs of positive integers $x$ and $y$ (with $x - y \neq 1$) such as to be able to identify the possible unique pair of non-trivial factors $a$ and $b$ of $n$ or the possible different pairs of non-trivial factors $a$ and $b$ of $n$, with $a > b$, starting from the value of $x = x_1$, as defined by the following equation:

$$x_1 = \frac{n - \left\lfloor \frac{n - 2 \cdot \lfloor \sqrt{n} \rfloor}{4} \right\rfloor \cdot 4 + 1}{2}, \tag{5. 1}$$

and proceeding, thereafter, always with the step size of 2 for the increase in the value of $x$, i.e. increasing the value of $x$ always by two units at a time (starting, therefore, from the 2nd iteration with the value of $x_2 = x_1 + 2$), until we are certain that pairs of positive integers $x$ and $y$ (with $x - y \neq 1$) can no longer occur such that we can identify pairs of non-trivial factors $a$ and $b$ of $n$ in addition to the pairs of non-trivial factors already identified, or, possibly, until we are certain of the primality of $n$:

It should be highlighted that, in the cases of $n = 7$, $n = 5$ and $n = 3$, the value of $x_1 = \frac{n - \left\lfloor \frac{n - 2 \cdot \lfloor \sqrt{n} \rfloor}{4} \right\rfloor \cdot 4 + 1}{2}$ coincides with the value of $x = \frac{n+1}{2}$ with which the pair of trivial factors of $n$ is identified.

Therefore, DFM-1 is still based on the repeated application of the equation $y = \sqrt{x^2 - n}$, but the difference from FFA-1 lies in the adoption of a step size of 2, instead of a step size of 1, for the increase in the value of $x$, starting from the 2nd iteration (with the value of $x_2 = x_1 + 2$) and until the identification of the pair of trivial factors of $n$, and in the determination of the initial value of $x$, i.e. $x_1$, which, if the conditions occur, may also coincide with the value of $\lceil \sqrt{n} \rceil$.

Below, we proceed with an example of the application of DFM-1, again with $n = 70399$ (as in the example of applying FFA-1 shown in Section 4 with Table 1), however highlighting that it is advisable, (as for the adoption of FFA-1), since it is more convenient, to stop at the identification of the first integer value of $y$ with $x - y \neq 1$, proceeding to check the factorization (or possible primality), by adopting DFM-1 again, of the first (or unique, but we cannot know a priori) identified non-trivial factors $a$ and $b$ of $n$.

Let $n = 70399$, which is the product of the primes 113, 89 and 7, be the composite to factorize. After checking that $n$ is not a perfect square, i.e. that the condition $x = \sqrt{70399}$ does not occur, and after having preliminarily subjected $n$ to a quick probabilistic test of primality, we proceed, by repeatedly applying of (2.7) (i.e. the equation $y = \sqrt{x^2 - n}$), with the search for the pairs of non-trivial factors $a$ and $b$ of $n$ until reaching the condition $x - y = 1$, starting from the value of $x = x_1 = \frac{n - \left\lfloor \frac{n - 2 \cdot \lfloor \sqrt{n} \rfloor}{4} \right\rfloor \cdot 4 + 1}{2} = 266$ and proceeding, thereafter, always with a step size of 2 for the increase in the value of $x$, with the results reported in Table 2 represented below, where the number of iterations is always denoted by $i$ (and with the value of $x_1$ which, in this case, coincides with the value of $\lceil\sqrt{n}\rceil$, while, if we had proceeded, for example, with the factorization of the composite $n$ = 70741, which is the product of the primes 109, 59 and 11, we would have had the value of $x_1 = \frac{n - \left\lfloor\frac{n - 2\cdot\lfloor\sqrt{n}\rfloor}{4}\right\rfloor\cdot 4 + 1}{2} = 267 > \lceil\sqrt{n}\rceil = 266$):

Table 2. Example of the application of DFM-1 with $n$ = 70399

| $i$ | $x$ | $y$ | $S$ | $a$ | $b$ |
|---|---|---|---|---|---|
| 1 | 266 | non-integer | 0 | ----- | ----- |
| 2 | 268 | non-integer | 0 | ----- | ----- |
| 3 | 270 | non-integer | 0 | ----- | ----- |
| 4 | 272 | non-integer | 0 | ----- | ----- |
| OMISSED | | | | | |
| 51 | 366 | non-integer | 0 | ----- | ----- |
| **52** | **368** | **255** | **0** | **623** | **113** |
| 53 | 370 | non-integer | 0 | ----- | ----- |
| 54 | 372 | non-integer | 0 | ----- | ----- |
| 55 | 374 | non-integer | 0 | ----- | ----- |
| OMISSED | | | | | |
| 87 | 438 | non-integer | 0 | ----- | ----- |
| **88** | **440** | **351** | **0** | **791** | **89** |
| 89 | 442 | non-integer | 0 | ----- | ----- |
| 90 | 444 | non-integer | 0 | ----- | ----- |
| 91 | 446 | non-integer | 0 | ----- | ----- |
| OMISSED | | | | | |
| 2383 | 5030 | non-integer | 0 | ----- | ----- |
| **2384** | **5032** | **5025** | **0** | **10057** | **7** |
| 2385 | 5034 | non-integer | 0 | ----- | ----- |
| 2386 | 5036 | non-integer | 0 | ----- | ----- |
| OMISSED | | | | | |
| 17467 | 35198 | non-integer | 0 | ----- | ----- |
| **17468** | **35200** | **35199** | **0** | **70399** | **1** |

As we can easily deduce from Table 2 and always denoting by $x_i$ the value of $x$ corresponding to a specific iteration considered, with DFM-1 we have:

$$i = \frac{x_i - x_1}{2} + 1 \tag{5.2}$$

for each iteration until the identification of the pair of trivial factors of $n$; therefore, with DFM-1 we always have:

$$i = \frac{\frac{a+b}{2} - x_1}{2} + 1 \tag{5.3}$$

both for the identification of any pair of non-trivial factors of $n$ and for the identification of the pair of trivial factors of $n$, and, in the latter case, we can also write the following equation:

$$i = \frac{\frac{n+1}{2} - x_1}{2} + 1 \quad (5.4)$$

since $a = n$ and $b = 1$, with the *value* of $i$ which, in the latter case, we can obviously know a priori.

That said, and as done for FFA-1, it is appropriate to highlight that, after the 1st iteration, the subsequent iterations can also be performed differently, i.e. without calculating the square of $x$ each time and without calculating the difference between $x^2$ and $n$ each time, thus reducing the computational cost (but always leaving the number of iterations to be performed unchanged). Indeed, starting from considering the known value (obviously integer) of $y^2$ corresponding to the 1st iteration and always denoting by $x_f$ the value of $x$ corresponding to the previous iteration and by $y_f^2$ and $y_r^2$, respectively, the value of $y^2 = x^2 - n$ corresponding to the previous iteration and the value of $y^2 = x^2 - n$ corresponding to the subsequent iteration, with DFM-1 we can also proceed by adopting the following equation:

$$y_r^2 = y_f^2 + 4 \cdot (x_f + 1), \quad (5.5)$$

for each iteration subsequent to the 1st iteration and until the identification of the pair of trivial factors of $n$ (starting, therefore, from the 2nd iteration with the value of $x_f = x_1 = \frac{n - \left\lfloor \frac{n - 2 \cdot \lfloor \sqrt{n} \rfloor}{4} \right\rfloor \cdot 4 + 1}{2}$ and, thereafter, increasing the value of $x_f$ always by two units at a time until the identification of the pair of trivial factors of $n$), and, then, checking whether the value of the square root of $y_r^2$, corresponding to the subsequent iteration, is an integer.

Furthermore, considering it appropriate to repeat what is described in the last part of Section 4 of this paper, it is possible to proceed to a further simplification, (which is already known, as indicated in reference [2]), reducing the number of square root operations to be performed and, therefore, further reducing the computational cost (but always leaving the number of iterations to be performed unchanged). Indeed, since, as we know, a perfect square can only end with specific digits (0, 1, 4, 5, 6, 9 in base 10), if the last digit of $y_1^2$ (where $y_1^2$ is the first value of $y_f^2$ possibly to be considered) or of $y_r^2$ (like the last digit of any $y^2 = x^2 - n$) is 2, 3, 7 or 8, we can eliminate the corresponding square root operation.

We conclude this section by highlighting again that, for this new factorization algorithm (DFM-1) (as for the substantial simplification, with or without the respective simplifications, which has not yet been published) and with or without the simplifications already presented in this section, the distinction (an easy and immediate distinction) of all odd natural numbers into two different categories (one of which excludes, with absolute certainty, the odd natural numbers which are perfect squares), although it always remains valid, is useless for the reasons that will be evident in a forthcoming paper.

## 6. Lemmas, Theorem, and First Properties of DFM-1

In what follows, we state the lemmas and the theorem of the new factorization algorithm (DFM-1).

**Lemma 1.** *Let n be any odd natural number other than a perfect square, and we denote by a and b, with* $a > b$*, both the factors of any possible pair of non-trivial factors of n (with the factors a and b which may not even be coprime) and the factors of the pair of trivial factors of n (therefore, with a and b both being odd positive integers, both when* $a \neq n$ *and* $b \neq 1$*, and when* $a = n$ *and* $b = 1$*). Replacing, for factorization purposes, the following equation, already given in Section 3, which defines the value of Euler's totient function* $\varphi(n)$ *of any n that is a semiprime (therefore, with* $p_1$ *and* $p_2$ *being distinct primes greater than 2, i.e. with* $p_1$ *and* $p_2$ *being different prime from each other and both odd):*

$$\varphi(n) = (p_1 - 1) \cdot (p_2 - 1) \tag{3.2b}$$

*with the following equation, already given in Section 3, which defines the value of the product* $(a - 1) \cdot (b - 1)$*, denoted by* $\varphi_s(n)$*, which is associated both with any possible pair of non-trivial factors a and b of n and with the pair of trivial factors a and b of n:*

$$\varphi_s(n) = (a - 1) \cdot (b - 1), \tag{3.3}$$

*and, therefore, by modifying, for factorization purposes, Euler's totient function* $\varphi(n)$ *of any n that is a semiprime, the two factors a and b of any possible pair of non-trivial factors of n and of the pair of trivial factors of n can be determined, using the corresponding value of* $\varphi_s(n)$*, via Fermat's Factorization Algorithm (FFA).*

*Proof.* From the following equation:

$$\varphi_s(n) = (a - 1) \cdot (b - 1), \tag{3.3}$$

that defines the value of $\varphi_s(n)$ associated both with any possible pair of non-trivial factors $a$ and $b$ of $n$ and with the pair of trivial factors $a$ and $b$ of $n$, we obtain the following equation which defines the sum $a + b$:

$$a + b = n - \varphi_s(n) + 1, \tag{6.1a}$$

i.e.:

$$\sum(a, b) = n - \varphi_s(n) + 1, \tag{6.1b}$$

since we have:

$$\begin{aligned} \varphi_s(n) &= (a - 1) \cdot (b - 1) \\ &= a \cdot b - a - b + 1 \\ &= n - (a + b) + 1 \end{aligned} \tag{6.2}$$

from which it follows that:

$$a + b = n - \varphi_s(n) + 1, \tag{6.1a}$$

i.e.:

$$\sum(a, b) = n - \varphi_s(n) + 1. \tag{6.1b}$$

From Fermat's Factorization Algorithm (FFA), we know that:

$$x = \frac{a+b}{2} \tag{2.3}$$

and, therefore, we have:

$$x = \frac{n - \varphi_s(n) + 1}{2} \tag{6.3}$$

and, consequently, Eq. (2.7) becomes:

$$y = \sqrt{\left(\frac{n - \varphi_s(n) + 1}{2}\right)^2 - n}, \tag{6.4}$$

with $x$ and $y$ both being positive integers (both with $x - y \neq 1$ and with $x - y = 1$), and, therefore, always using Fermat's Factorization Algorithm (FFA), we can determine, as indicated below and as we already know from Section 2, the two factors $a$ and $b$ of any possible pair of non-trivial factors of $n$ and of the pair of trivial factors of $n$:

$$a = x + y \tag{2.9}$$

and

$$b = x - y, \tag{2.10}$$

so that, by adopting (6.3) and (6.4), we have the two factors $a$ and $b$ of any possible pair of non-trivial factors of $n$ and of the pair of trivial factors of $n$ as determined below, using the corresponding value of $\varphi_s(n)$, via Fermat's Factorization Algorithm (FFA):

$$a = \frac{n - \varphi_s(n) + 1}{2} + \sqrt{\left(\frac{n - \varphi_s(n) + 1}{2}\right)^2 - n} \tag{6.5}$$

and

$$b = \frac{n - \varphi_s(n) + 1}{2} - \sqrt{\left(\frac{n - \varphi_s(n) + 1}{2}\right)^2 - n}, \tag{6.6}$$

and, hence, Lemma 1 is proved. □

**Lemma 2.** *Let n be any odd natural number other than a perfect square, if n is a semi-prime or non-semiprime composite and denoting by a and b, with a > b, the factors, prime or composite (which may not even be coprime), of any pair of non-trivial factors of n (therefore, with a and b both being odd positive integers, where a ≠ n and b ≠ 1), the value of*

$$\varphi_s(n) = (a - 1) \cdot (b - 1), \tag{3.3}$$

*which is associated with any pair of non-trivial factors a and b of of n, is an even positive integer that is always strictly less than the following value of* $n - 2 \cdot \lfloor\sqrt{n}\rfloor$, *denoted by* $\varphi_{s1_p}(n)$, *which is always an odd positive integer:*

$$\varphi_{s1_p}(n) = n - 2 \cdot \lfloor\sqrt{n}\rfloor, \tag{6.7}$$

*and, therefore, we always have:*

$$\varphi_s(n) < \varphi_{s1_p}(n), \tag{6.8a}$$

*i.e.:*

$$(a - 1) \cdot (b - 1) < n - 2 \cdot \lfloor\sqrt{n}\rfloor. \tag{6.8b}$$

*Proof.* First of all, it is worth pointing out that $\varphi_{s1_p}(n) = n - 2 \cdot \lfloor\sqrt{n}\rfloor$ is always an odd positive integer. Indeed, since $n$ is an odd natural number and the product $2 \cdot \lfloor\sqrt{n}\rfloor$ is obviously always an even positive integer (whether $\lfloor\sqrt{n}\rfloor$ is an even positive integer or an odd positive integer), the difference $n - 2 \cdot \lfloor\sqrt{n}\rfloor$, obviously, is always an odd positive integer; furthermore, since a and b are both odd positive integers other than 1, $(a - 1)$ and $(b - 1)$ are both even positive integers and, consequently, the product $(a - 1) \cdot (b - 1)$ is always an even positive integer. This, therefore, rules out a priori the condition $\varphi_s(n) = \varphi_{s1_p}(n)$.

That said, and since we have, as we already know from Lemma 1:

| | |
|---|---|
| $\varphi_s(n)$ | $= (a - 1) \cdot (b - 1)$ |
| | $= a \cdot b - a - b + 1$ |
| | $= n - (a + b) + 1,$ |

(6.2)

Relation (6.8*b*) becomes:

$$n - (a + b) + 1 < n - 2 \cdot \lfloor\sqrt{n}\rfloor, \tag{6.8c}$$

which, simplifying, becomes:

$$a + b > 2 \cdot \lfloor\sqrt{n}\rfloor + 1. \tag{6.8d}$$

Therefore, we must prove that, with $a$ and $b$ both being odd positive integers and with $a \neq b$, Relation (6.8*d*) always holds.

First of all, it is easy to verify that, for (6.8*d*) to always hold, as well as (6.8*b*), both $a$ and $b$ must be odd positive integers (or both even positive integers, if, by hypothesis, $n$ were not exclusively any odd natural number other than a perfect square) since, if, by hypothesis, $n$ were not exclusively any odd natural number other than a perfect square and $a$ and $b$ did not have the same parity (i.e. one were even and the other were odd), and $a$ and $b$ were very close to each other (relative to the order of magnitude of their product $a \cdot b$), Relations (6.8*b*) and (6.8*d*) would become, respectively:

$$(a-1) \cdot (b-1) = a \cdot b - 2 \cdot \lfloor \sqrt{a \cdot b} \rfloor$$

and

$$a + b = 2 \cdot \lfloor \sqrt{a \cdot b} \rfloor + 1,$$

with the product $a \cdot b$ which is, unlike $n$, an even natural number. So, in the case in question, $a$ and $b$ are both odd positive integers since $n$ is any odd natural number other than a perfect square.

That said, we know that, with $a$ and $b$ both being odd positive integers and with $a \neq b$, the AM-GM (arithmetic mean-geometric mean) inequality states that the arithmetic mean is strictly greater than the geometric mean, that is:

$$\frac{a+b}{2} > \sqrt{a \cdot b} = \sqrt{n}, \tag{6.9a}$$

i.e.:

$$a + b > 2 \cdot \sqrt{a \cdot b} = 2 \cdot \sqrt{n}, \tag{6.9b}$$

having excluded a priori (since $n$ is any odd natural number other than a perfect square) the relation $\frac{a+b}{2} = \sqrt{a \cdot b} = \sqrt{n}$, i.e. the relation $a + b = 2 \cdot \sqrt{a \cdot b} = 2 \cdot \sqrt{n}$, which occurs with $a = b$.

Since $\sqrt{a \cdot b}$ is greater than $\lfloor \sqrt{a \cdot b} \rfloor$ being $\sqrt{a \cdot b}$ non-integer, and since, by the AM-GM inequality, $\frac{a+b}{2}$ is greater than $\sqrt{a \cdot b}$, we have:

$$\frac{a+b}{2} > \sqrt{a \cdot b} > \lfloor \sqrt{a \cdot b} \rfloor$$

and, therefore, we have:

$$\frac{a+b}{2} > \lfloor \sqrt{a \cdot b} \rfloor, \tag{6.10a}$$

i.e.:

$$a + b > 2 \cdot \lfloor\sqrt{a \cdot b}\rfloor. \tag{6.10b}$$

Since, trivially, $\frac{a+b}{2}$ is greater than $\lfloor\sqrt{a \cdot b}\rfloor$ by at least 1, it follows, trivially, that $a + b$ is greater than $2 \cdot \lfloor\sqrt{a \cdot b}\rfloor$ by at least 2, and, therefore, we have the following relation:

$$a + b \geq 2 \cdot \lfloor\sqrt{a \cdot b}\rfloor + 2 = 2 \cdot \lfloor\sqrt{n}\rfloor + 2, \tag{6.11}$$

which, since we are dealing with integers, can also be written as follows, thus holding (6.8*d*):

$$a + b > 2 \cdot \lfloor\sqrt{n}\rfloor + 1, \tag{6.8d}$$

and, therefore, the following (6.8*b*) holds:

$$(a - 1) \cdot (b - 1) < n - 2 \cdot \lfloor\sqrt{n}\rfloor, \tag{6.8b}$$

and, hence, Lemma 2 is proved. □

**Lemma 3.** *Let n be any odd natural number other than a perfect square, if n is a semiprime or non-semiprime composite and denoting by a and b, with a > b, the factors, prime or composite (which may not even be coprime), of any pair of non-trivial factors of n (therefore, with a and b both being odd positive integers, where* $a \neq n$ *and* $b \neq 1$*), the value of*

$$\varphi_s(n) = (a - 1) \cdot (b - 1), \tag{3.3}$$

*which is associated with any pair of non-trivial factors a and b of n, is an even positive integer that is certainly always a multiple of* 4 *and, hence, certainly always divisible by* 4.

*Proof.* Since $a$ and $b$ are both odd positive integers other than 1, $(a - 1)$ and $(b - 1)$ are both even positive integers and, therefore, $(a - 1)$ and $(b - 1)$ are both certainly formed, among other possible odd prime factors, by the factor 2 with exponent at least equal to 1 and, consequently, the product $(a - 1) \cdot (b - 1)$ is an even positive integer which is certainly formed, among other possible odd prime factors, by the factor 2 with exponent at least equal to 2; therefore, the value of the product $(a - 1) \cdot (b - 1)$, which is the value of $\varphi_s(n)$ associated with any pair of non-trivial factors $a$ and $b$ of any $n$ that is a semiprime or non-semiprime composite, is an even positive integer that is certainly always a multiple of 4 and, hence, certainly always divisible by 4.

Hence, Lemma 3 is proved. □

**Lemma 4.** *Let n be any odd natural number other than a perfect square, and we denote by a and b, with a > b, the factors of the pair of trivial factors of any n (therefore, with a and b both being odd positive integers, where* $a = n$ *and* $b = 1$*). The value of*

$$\varphi_s(n) = (a - 1) \cdot (b - 1), \tag{3.3}$$

*which is associated with the pair of trivial factors a and b of n, is equal to* 0*, which, as such, is certainly also a multiple of* 4*.*

*Proof.* Since $a = n$ and $b = 1$, we have the following value of $\varphi_s(n) = (a-1) \cdot (b-1)$ that is associated with the pair of trivial factors $a$ and $b$ of $n$:

$$\begin{aligned}\varphi_s(n) &= (a-1) \cdot (b-1) \\ &= (n-1) \cdot (1-1) \\ &= 0,\end{aligned} \tag{6.12}$$

with 0 which, as we know, is certainly also a multiple of 4.

Hence, Lemma 4 is proved. □

**Lemma 5.** *Let n be any odd natural number other than a perfect square, and we denote by a and b, with a > b, the factors of the pair of trivial factors of any n (therefore, with a and b both being odd positive integers, where a = n and b =* 1*). If n is a prime, only the value of* $\varphi_s(n) = (a-1) \cdot (b-1) = 0$ *occurs, which is associated with the pair of trivial factors a and b of n and which, as such, is certainly also a multiple of* 4*.*

*Proof.* Since $a = n$ and $b = 1$ are the only integer factors of $n$ if $n$ is a prime, only the value of $\varphi_s(n) = (a-1) \cdot (b-1) = 0$ can occur, i.e. the value of $\varphi_s(n) = (a-1) \cdot (b-1)$ that is associated with the pair of trivial factors $a$ and $b$ of $n$ and that, as we know, is certainly also a multiple of 4, with the proof that is identical to that of Lemma 4, and, hence, Lemma 5 is proved. □

**Lemma 6.** *Let n be any odd natural number other than a perfect square, and we denote by a and b, with a > b, both the factors of any possible pair of non-trivial factors of n (with the factors a and b which may not even be coprime) and the factors of the pair of trivial factors of n (therefore, with a and b both being odd positive integers, both when a ≠ n and b ≠* 1*, and when a = n and b =* 1*). After verifying, by Lemmas 1 and 2, the value of y associated with the following computed initial value of* $\varphi_s(n)$*, which we denote by* $\varphi_{s1}(n)$ *and which is, therefore, the largest possible value of* $\varphi_s(n)$*, i.e. the value of the* largest *multiple of* 4 *(which, by Lemmas 3 and 4, is an even positive integer or equal to* 0*) less than* $\varphi_{s1_p}(n) = n - 2 \cdot \lfloor\sqrt{n}\rfloor$*, which, as indicated in Lemma 2, is always an odd positive integer:*

$$\varphi_{s1}(n) = \left\lfloor \frac{\varphi_{s1_p}(n)}{4} \right\rfloor \cdot 4 = \left\lfloor \frac{n - 2 \cdot \lfloor\sqrt{n}\rfloor}{4} \right\rfloor \cdot 4 \tag{6.13}$$

*and, therefore, considering (6.3) and (6.4), with the following initial values of x and y, which we denote, respectively, by* $x_1$ *and* $y_1$*, and with the value of* $x_1$ *already known from Section 5:*

$$x_1 = \frac{n - \varphi_{s1}(n) + 1}{2} \tag{6.14}$$

*and*

$$y_1 = \sqrt{\left(\frac{n - \varphi_{s1}(n) + 1}{2}\right)^2 - n}, \tag{6.15}$$

*provided that* $\varphi_{s1}(n) \neq 0$ *and always via Fermat's Factorization Algorithm (FFA) (i.e. by repeatedly applying the equation* $y = \sqrt{x^2 - n}$*), in order to identify any possible integer value of y with* $x - y \neq 1$ *(or any possible further integer value of y with* $x - y \neq 1$*, if, with the value of* $\varphi_{s1}(n)$*, an integer value of y with* $x - y \neq 1$ *has been identified, which could also be the only integer value of y with* $x - y \neq 1$*)—i.e. any possible pair of non-trivial factors a and b of n (or any possible further pair of non-trivial factors a and b of n, if, with the value of* $\varphi_{s1}(n)$*, a pair of non-trivial factors a and b of n has been identified, which could also be the only pair of non-trivial factors a and b of n)—and the integer value of y with* $x - y = 1$ *(i.e. the pair of trivial factors a and b of n), or in order to determine the primality of n, the value of* $\varphi_{s1}(n) = \left\lfloor \frac{n - 2 \cdot \lfloor\sqrt{n}\rfloor}{4} \right\rfloor \cdot 4$ *is decreased, by Lemmas 3 and 4, by four units at a time in the following equation:*

$$y = \sqrt{\left(\frac{n - \varphi_{s1}(n) + 1}{2}\right)^2 - n}, \tag{6.16}$$

*or, to simplify, the value of* $x_1 = \frac{n - \varphi_{s1}(n) + 1}{2}$ *is increased by two units at a time in the following equation:*

$$y = \sqrt{x_1^2 - n}. \tag{6.17}$$

*Proof.* Let $\varphi_{s1}(n) \neq 0$ and we denote by $\varphi_{sX}(n)$ the value of $\varphi_s(n)$ associated with any possible pair of non-trivial factors *a* and *b* of *n* and with the pair of trivial factors *a* and *b* of *n* (hence, both with $\varphi_{sX}(n) \neq 0$ and with $\varphi_{sX}(n) = 0$). Since the value of $\varphi_{s1}(n)$, with $\varphi_{s1}(n) \neq 0$, is an even positive integer that is a multiple of 4, and, by Lemmas 3 and 4, the value of $\varphi_{sX}(n)$ is certainly a multiple of 4 (which, again by Lemmas 3 and 4, is an even positive integer or equal to 0), it is immediate that, with $\varphi_{s1}(n) > \varphi_{sX}(n)$, their difference is an even positive integer that is a multiple of 4, and that, starting from the value of $\varphi_{s1}(n)$ (which could also coincide with the value of $\varphi_s(n)$ associated with the first or unique pair of non-trivial factors *a* and *b* of *n*), we can reach the value of $\varphi_{sX}(n)$ by decreasing the value of $\varphi_{s1}(n)$ by four units at a time until the entire difference $\varphi_{s1}(n) - \varphi_{sX}(n)$ is covered.

At this point, we must make another consideration: denoting by $x_2$ the value of *x* corresponding to the 2nd iteration, decreasing the value of $\varphi_{s1}(n) = \left\lfloor \frac{n - 2 \cdot \lfloor\sqrt{n}\rfloor}{4} \right\rfloor \cdot 4$ by four units means increasing the value of $x_1 = \frac{n - \varphi_{s1}(n) + 1}{2}$ by two units, as is evident in the following simplification:

$$
\begin{aligned}
x_2 &= \frac{n - (\varphi_{s1}(n) - 4) + 1}{2} \\
&= \frac{n - \varphi_{s1}(n) + 4 + 1}{2} \\
&= \frac{n - \varphi_{s1}(n) + 1}{2} + \frac{4}{2} \\
&= x_1 + 2,
\end{aligned} \tag{6.18}
$$

with the following corresponding value of $y_2$:

$$y_2 = \sqrt{x_2{}^2 - n}, \tag{6.19a}$$

i.e.:

$$y_2 = \sqrt{(x_1 + 2)^2 - n} = \sqrt{\left(\frac{n - \varphi_{s1}(n) + 1}{2} + 2\right)^2 - n}, \tag{6.19b}$$

and, therefore, it is immediate that the values of $x$ increase by two units with each subsequent iteration.

Hence, Lemma 6 is proved. □

**Lemma 7.** *Let n be any odd natural number other than a perfect square and different from* 3, 5, *and* 7, *and we denote by a and b, with* $a > b$, *the factors, both being odd positive integers or both being positive non-integers, of any pair of factors of n (therefore, with* $a \neq n$ *and* $b \neq 1$, *where a and b are both odd positive integers and may not even be coprime; with a and b both being odd positive integers, where* $a = n$ *and* $b = 1$; *and also with a and b both being positive non-integers such that* $a \cdot b = n$*). Provided that* $\varphi_{s1}(n) \neq 0$, *as we increase the value of* $x_1 = \frac{n - \varphi_{s1}(n) + 1}{2}$ *by two units at a time in the following equation according to Lemma 6:*

$$y = \sqrt{x_1{}^2 - n}, \tag{6.17}$$

*the value of a increases more and more and the value of b decreases more and more, i.e. with a and b moving further and further away from each other, until the identification of the pair of trivial factors of n, with the value of* $\varphi_s(n) = (a - 1) \cdot (b - 1) = 0$ *(by Lemma 4), with the largest possible integer value of a and the smallest possible integer value of b (i.e. with a and b, both being odd positive integers, as far apart as possible), being* $a = n$ *and* $b = 1$.

*Proof.* First of all, it should be noted that, in this Lemma 7, $n$ is any odd natural number other than both a perfect square and 3, 5, and 7; this is due to the fact that, with $n$ being any odd natural number other than a perfect square, we have obviously excluded the cases with the condition $\varphi_s(n) = \varphi_{s1}(n) = 0$, i.e. the condition $x_1 = \frac{n - \varphi_{s1}(n) + 1}{2} = \frac{n + 1}{2}$, and for which, therefore, no iterations subsequent to the 1st iteration occur.

That said, and to simplify the proof of this Lemma 7, we assume that $n$ is a semiprime or non-semiprime composite and we denote by $a_1$ and $b_1$, with $a_1 > b_1$, the factors of the first or unique pair of non-trivial factors of $n$; and we also assume that the value of $\varphi_{s1}(n)$ coincides with the value of $\varphi_s(n) = (a_1 - 1) \cdot (b_1 - 1)$ associated with the first or unique pair of non-trivial factors of $n$. By Lemma 1 and according to what is established in Lemma 6 regarding the value of $\varphi_{s1}(n)$, we have:

$$a_1 = \frac{n - \varphi_{s1}(n) + 1}{2} + \sqrt{\left(\frac{n - \varphi_{s1}(n) + 1}{2}\right)^2 - n}, \tag{6.20a}$$

i.e.:

$$a_1 = x_1 + y_1 = x_1 + \sqrt{x_1{}^2 - n}, \tag{6.20b}$$

with $y_1$ which is obviously an integer since $a_1$ and $b_1$ are both integers, as we assumed for this proof that $a_1$ and $b_1$ are the factors of the first or unique pair of non-trivial factors of $n$.

At this point, by Lemma 6 and denoting by $a_2$and $b_2$, with $a_2 > b_2$, the factors—both being odd positive integers or both being positive non-integers—of the pair of factors of $n$ determined by the 2nd iteration, i.e. determined by the first increment in the value of $x_1 = \frac{n - \varphi_{s1}(n) + 1}{2}$ by two units, we have:

$$a_2 = x_2 + y_2 = (x_1 + 2) + \sqrt{(x_1 + 2)^2 - n}, \tag{6.21}$$

from which it clearly follows that $a_2$ is greater than $a_1$, since it is clearly evident that:

$$x_2 = (x_1 + 2) > x_1 \tag{6.22}$$

and

$$y_2 = \sqrt{(x_1 + 2)^2 - n} > y_1 = \sqrt{x_1{}^2 - n}, \tag{6.23}$$

and, therefore, we have:

$$a_2 > a_1 \tag{6.24a}$$

since it is clearly evident that:

$$\begin{aligned} a_2 &= x_2 + y_2 \\ &= x_2 + \sqrt{x_2{}^2 - n} \\ &= (x_1 + 2) + \sqrt{(x_1 + 2)^2 - n} \\ &> a_1 \\ &= x_1 + y_1 \\ &= x_1 + \sqrt{x_1{}^2 - n}. \end{aligned} \tag{6.24b}$$

At this point, since we know $a_1 \cdot b_1 = a_2 \cdot b_2 = n$ and we have shown that $a_2 > a_1$, we can unequivocally conclude that:

$$b_2 < b_1 \tag{6.25}$$

since we have:

$$b_2 = \frac{n}{a_2} < b_1 = \frac{n}{a_1}, \tag{6.26}$$

with $b_2$ being necessarily and logically less than $b_1$ since the two fractions given above in (6.26) have the same numerator $n$, and, therefore, $b_2$ must be less than $b_1$ since the denominator $a_2$ of the fraction $\frac{n}{a_2}$ is greater than the denominator$a_1$ of the fraction $\frac{n}{a_1}$, and with $a_2$ and $b_2$ for which the following equation always holds (where $i$ is the number of iterations, inclusive of the 1st iteration, which in this case is $i = 2$), with the sum $\sum(a_2, b_2)$ which is always an even positive integer:

$$\sum(a_2, b_2) = n - \varphi_{s1}(n) + (i - 1) \cdot 4 + 1, \tag{6.27$a$}$$

i.e.:

$$\sum(a_2, b_2) = \sum(a_1, b_1) + (i - 1) \cdot 4, \tag{6.27$b$}$$

and also:

$$\sum(a_2, b_2) = n - (a_2 - 1) \cdot (b_2 - 1) + 1, \tag{6.27$c$}$$

even though it may happen that, with the first increment in the value of $x_1 = \frac{n - \varphi_{s1}(n) + 1}{2}$ by two units, $a_2$ and $b_2$, unlike $a_1$ and $b_1$ (which in this proof are both integers, given the assumption that they are the factors of the first or unique pair of non-trivial factors of $n$), are not integers, since $y_2$, unlike $y_1$, turns out to be non-integer (but, obviously, with the sum $\sum(a_2, b_2)$ which is always an even positive integer, being—see the right-hand sides of (6.27a), (6.27b), and (6.27c)—$n - \varphi_{s1}(n) + (i - 1) \cdot 4 + 1$, i.e. $\sum(a_1, b_1) + (i - 1) \cdot 4$ and also $n - (a_2 - 1) \cdot (b_2 - 1) + 1$, always an even positive integer, with the product $(a_2 - 1) \cdot (b_2 - 1)$ which, like the product $(a_1 - 1) \cdot (b_1 - 1)$ and as already pointed out in Lemma 2, is, therefore, always an even positive integer). This means that with DFM-1, given that the value of $x$ is incremented by two units at a time, it can happen, if $n$ is a non-semiprime composite (as, for example, in the cases of $n = 3^3 \cdot 5 \cdot 7 = 945$ and $n = 127 \cdot 129 \cdot 131 \cdot 133 = 285441009$), that $a_2$ and $b_2$, like $a_1$ and $b_1$, are integers (with $a_2 \neq n$ and $b_2 \neq 1$, where $a_2$ and $b_2$ are both odd positive integers and may not even be coprime as in the case of the second pair of non-trivial factors of $n = 945$, where $a_2 = 3^2 \cdot 5 = 45$ and $b_2 = 3 \cdot 7 = 21$) since $y_2$, like $y_1$, turns out to be an integer with the condition $x_2 - y_2 \neq 1$, and, therefore, it may happen that the first and second pairs of non-trivial factors of $n$ are identified, respectively, at the 1st and 2nd iterations; however, this circumstance can occur only if the first pair of non-trivial factors of $n$ is identified at the 1st iteration, i.e. when the value of $\varphi_{s1}(n) = \left\lfloor \frac{n - 2 \cdot \lfloor\sqrt{n}\rfloor}{4} \right\rfloor \cdot 4$ coincides with the value of $\varphi_s(n) = (a_1 - 1) \cdot (b_1 - 1)$ associated with the first pair of non-trivial factors of $n$, and, consequently, with $a_1$ and $b_1$ being determined by the value of $x_1 = \frac{n - \varphi_{s1}(n) + 1}{2}$ and the corresponding value of $y_1$; nevertheless, this circumstance does not occur in all cases where the first pair of non-trivial factors of $n$—with $n$ thus being a non-semiprime composite—is identified at the 1st iteration. While, with FFA-1, given that the value of $x$ is incremented by one unit at a time, it can never happen—for the reasons that are trivially already evident in this paper—that, with $n$ which is a non-semiprime composite, the first and second pairs of non-trivial factors of $n$ are identified, respectively, at the 1st and 2nd iterations.

Obviously, this proof (which, clearly, also holds for FFA-1, with the value of $x$ increasing by one unit at a time) also holds both for the case where $a_1$ and $b_1$ were not integers and for further increases in the value of $x$ by two units at a time [with the corresponding values of $\sum(a, b)$ which are, obviously, always even positive integers, being, for any iteration, $n - \varphi_{s1}(n) + (i - 1) \cdot 4 + 1$, i.e. $\sum(a_1, b_1) + (i - 1) \cdot 4$ and also $n - (a - 1) \cdot (b - 1) + 1$, always an even positive integer, with the product $(a - 1) \cdot (b - 1)$ that is, therefore, always an even positive integer for any iteration], confirming that, as we increase the value of $x_1 = \frac{n - \varphi_{s1}(n) + 1}{2}$ by two units at a time in (6.17) (i.e. the equation $y = \sqrt{x_1{}^2 - n}$), the value of $a$ (even if non-integer) increases more and more and the value of $b$ (even if non-integer) decreases more and more, i.e. with $a$ and $b$ moving further and further away from each other, until the identification of the pair of trivial factors of $n$, with the value of $\varphi_s(n) = (a - 1) \cdot (b - 1) = 0$, with the largest possible integer value of $a$ and the smallest possible integer value of $b$ (i.e. with $a$ and $b$, both being odd positive integers, as far apart as possible), being $a = n$ and $b = 1$ [obviously, as we increase the value of $x_1 = \frac{n - \varphi_{s1}(n) + 1}{2}$ by two units at a time in the equation $y = \sqrt{x_1{}^2 - n}$, the value of the product $(a - 1) \cdot (b - 1)$ decreases more and more until the identification of the pair of trivial factors $a$ and $b$ of $n$ with the value of $\varphi_s(n) = (a - 1) \cdot (b - 1) = 0$, while the values of $\sum(a, b)$ increases more and more until it reaches the value $\sum(a, b) = n + 1$, which is associated with the pair of trivial factors $a$ and $b$ of $n$].

Hence, Lemma 7 is proved. □

**Theorem.** *Let n be any odd natural number other than a perfect square, and we denote by a and b, with a > b, both the factors of any possible pair of non-trivial factors of n (with the factors a and b which may not even be coprime) and the factors of the pair of trivial factors of n (therefore, with a and b both being odd positive integers, both when $a \neq n$ and $b \neq 1$, and when $a = n$ and $b = 1$). By Lemmas 1-7, by modifying, for factorization purposes, Euler's totient function $\varphi(n)$ of any n that is a semiprime, and using Fermat's Factorization Algorithm (FFA)—i.e. by repeatedly applying the equation $y = \sqrt{x^2 - n}$—provided that $x_1 \neq \frac{n+1}{2}$ (being $x_1 = \frac{n+1}{2}$ the value of x associated with the pair of trivial factors of n), any possible pair of non-trivial factors a and b of n can be identified by looking for the possible corresponding integer value of y (obviously, with $x - y \neq 1$), starting from the value of $x = x_1 = \frac{n - \left\lfloor \frac{n - 2 \cdot \lfloor\sqrt{n}\rfloor}{4} \right\rfloor \cdot 4 + 1}{2}$ and proceeding, thereafter, always with the step size of 2, i.e. increasing the value of x always by two units at a time (starting, therefore, from the 2nd iteration with the value of $x_2 = x_1 + 2$), until we are certain that pairs of positive integers x and y (with $x - y \neq 1$) can no longer occur such that we can* *identify pairs of non-trivial factors a and b of n in addition to the pairs of non-trivial factors already identified, or, possibly, until we are certain of the primality of n.*

*If n is a semiprime composite, the only integer value of y found (with* $x - y \neq 1$*) is the one that determines, with the corresponding value of x, the unique pair of non-trivial factors a and b of n, with the non-trivial value of a which is smaller than the trivial value of a and the non-trivial value of b which is larger than the trivial value of b (i.e. with non-trivial a and b closest to each other compared to what happens with the pair of trivial factors a and b of n).*

*While, if n is a non-semiprime composite, the first integer value of y found (with* $x - y \neq 1$*) is the one that determines, with the corresponding value of x, the first pair of non-trivial factors a and b of n, with the smallest non-trivial value of a and the largest non-trivial value of b (i.e. with non-trivial a and b closest to each other), and the last integer value of y found (with* $x - y \neq 1$*) is the one that determines, with the corresponding value of x, the last pair of non-trivial factors a and b of n, with the largest non-trivial value of a and the smallest non-trivial value of b (i.e. with non-trivial a and b farthest from each other). So, if n is a non-semiprime composite, as we identify the pairs of non-trivial factors a and b of n, the value of a increases more and more and the value of b decreases more and more, i.e. with a and b moving further and further away from each other, until the identification of the pair of trivial factors a and b of n, with the largest possible value of a and the smallest possible value of b (i.e. with a and b as far apart as possible, being* $a = n$ and $b = 1$*), which is, obviously, the only pair of (integer) factors a and b of n that occurs if n is a prime.*

*Proof.* The proof is deferred to Lemmas 1 through 7. □

Therefore, after checking that the condition $x = \sqrt{n}$ does not occur and after having preliminarily subjected $n$ to a quick probabilistic test of primality, we will proceed, by repeatedly applying of (2.7) (i.e. the equation $y = \sqrt{x^2 - n}$)—and possibly using only the equation $y^2 = x^2 - n$—and after checking the value of $y$ (always using the equation $y = \sqrt{x^2 - n}$ and possibly using only the equation $y^2 = x^2 - n$) for $x = x_1 = \frac{n - \left\lfloor \frac{n - 2 \cdot \lfloor \sqrt{n} \rfloor}{4} \right\rfloor \cdot 4 + 1}{2}$, always with the step size of 2, i.e. we will proceed by increasing, always by two units at a time, the value of $x$ (starting, therefore, from the 2nd iteration with the value of $x_2 = x_1 + 2 = \frac{n - \left\lfloor \frac{n - 2 \cdot \lfloor \sqrt{n} \rfloor}{4} \right\rfloor \cdot 4 + 1}{2} + 2$) to find possible integer values of $y$ with $x - y \neq 1$ and until we are certain that an integer values of $y$, with $x - y \neq 1$, can no longer occur such as to be able to identify pairs of non-trivial factors $a$ and $b$ of $n$ in addition to the pairs of non-trivial factors already identified, or, possibly, until we are certain of the primality of $n$.

Or (resorting to the simplification described in the last part of Section 5 of this paper), after checking that the condition $x = \sqrt{n}$ does not occur and after having preliminarily subjected $n$ to a quick probabilistic test of primality, we can also proceed, after checking the value of $y$ (using the equation $y = \sqrt{x^2 - n}$ and possibly using only the equation $y^2 = x^2 - n$) for $x = x_1 = \frac{n - \left\lfloor \frac{n - 2 \cdot \lfloor \sqrt{n} \rfloor}{4} \right\rfloor \cdot 4 + 1}{2}$, by adopting the equation $y_r^{\,2} = y_f^{\,2} + 4 \cdot (x_f + 1)$ as given in Section 5 with (5.5), for each iteration subsequent to the 1st iteration (starting, therefore, from the 2nd iteration with the value of $x_f = x_1 = \frac{n - \left\lfloor \frac{n - 2 \cdot \lfloor \sqrt{n} \rfloor}{4} \right\rfloor \cdot 4 + 1}{2}$ and, thereafter, increasing the value of $x_f$ always by two units at a time), to find possible integer values of $y_r$ with $x_f + 2 - y_r \neq 1$, with the values of $x_f$ and $y_r$ that appear in the same equation $y_r^2 = y_f^2 + 4 \cdot (x_f + 1)$, and until we are certain that an integer values of $y_r$, with $x_f + 2 - y_r \neq 1$, can no longer occur such as to be able to identify pairs of non-trivial factors $a$ and $b$ of $n$ in addition to the pairs of non-trivial factors already identified, or, possibly, until we are certain of the primality of $n$.

However, it is worth noting again that it is advisable (as for the adoption of FFA-1, and with or without the simplifications described in the last part of Section 5 of this paper), since it is more convenient, to stop at the identification of the first (which could also be the only one) integer value of $y$ with $x - y \neq 1$ (or integer value of $y_r$ with $x_f + 2 - y_r \neq 1$), provided that there exists an integer value of $y$ with $x - y \neq 1$ (or an integer value of $y_r$ with $x_f + 2 - y_r \neq 1$), proceeding to check the factorization (or possible primality), by adopting DFM-1 again (with or without the simplifications described in the last part of Section 5 of this paper), of the first (or unique, but we cannot know a priori) identified non-trivial factors $a$ and $b$ of $n$ (obviously, after checking that $a$ and $b$ are not perfect squares and after having preliminarily subjected $a$ and $b$ to a quick probabilistic test of primality).

**Note.** We will always stop at identifying the possible first or unique integer value of $y$ with $x - y \neq 1$ (or the possible first or unique integer value of $y_r$ with $x_f + 2 - y_r \neq 1$). Therefore, we will always focus primarily on identifying the possible first or unique pair of non-trivial factors of $n$.

Therefore, having found, possibly, the first or unique pair of positive integers $x$ and $y$ with $x - y \neq 1$ (or of positive integers $x_f + 2$ and $y_r$ with $x_f + 2 - y_r \neq 1$), we will proceed to determine the non-trivial factors $a$ and $b$ of $n$ using the following two known equations (highlighting that, with n for which more than one pair of non-trivial factors occurs, the first identified pair of non-trivial factors $a$ and $b$ of $n$ is the one with the smallest non-trivial value of $a$ and the largest non-trivial value of $b$, i.e. with non-trivial $a$ and $b$ closest to each other; while, the unique pair of non-trivial factors $a$ and $b$ of $n$ occurs, obviously, with $n$ being a semiprime or with $n$ being a cube of prime):

$$a = x + y \tag{2.9}$$

and

$$b = x - y, \tag{2.10}$$

or, resorting to the simplification described in the last part of Section 5 of this paper, by adopting the equation $y_r^2 = y_f^2 + 4 \cdot (x_f + 1)$ as given in Section 5 with (5.5), using the following two equations:

$$a = x_f + 2 + y_r \tag{6.28}$$

and

$$b = x_f + 2 - y_r, \tag{6.29}$$

and, subsequently, we will proceed to check the factorization (or possible primality) of these identified non-trivial factors $a$ and $b$ of $n$, by adopting DFM-1 again, with or without the simplifications described in the last part of Section 5 of this paper (obviously, after checking that $a$ and $b$ are not perfect squares and after having preliminarily subjected $a$ and $b$ to a quick probabilistic test of primality).

However, it is also advisable (as for the adoption of di FFA-1) to resort to the Trial Division Algorithm (TDA) [5, pp. 133-134], before proceeding with DFM-1 (with or without the simplifications described in the last part of Section 5 of this paper), attempting to divide $n$ by a relatively small number (relative to the order of magnitude of $n$) of successive primes, starting from the number 3, contained in $\sqrt{n}$ and possibly finding prime factors of small magnitude so as to possibly have to proceed to check the factorization (or possible primality) of an $n$ smaller than the one initially considered. And, as will be seen later, proceeding a priori in this sense is more convenient in cases of the primality of $n$ in order to reach, with the adoption of DFM-1 (with or without the simplifications described in the last part of Section 5 of this paper), the limit number of iterations, to be specified later, with a smaller number of iterations.

So, going into details of the new factorization algorithm (DFM-1), and also considering it appropriate to repeat what has already been described in the Lemmas and in the Theorem that have been stated in this section, we start by considering the following value of $\varphi(n)$, which is the value of Euler's totient function [4] of any $n$ that is a semiprime, which from now on we replace, for factorization purposes, with the value of $\varphi_s(n)$ to whose formation contribute the possible factors $a$ and $b$, prime or composite (which may not even be coprime) and with $a > b$, which form the possible first or unique pair of non-trivial factors $a$ and $b$ of $n$ and which are, therefore, in any case always considered, for the purposes of applying DFM-1, as primes even if they are not (this holds both for any further possible pair of non-trivial factors $a$ and $b$ of $n$ and for the pair of trivial factors $a$ and $b$ of $n$, always with $a > b$), thus putting a particular light on Euler's totient function, a particular light which consists in the modification, for factorization purposes, of Euler's totient function $\varphi(n)$ of any $n$ that is a semiprime:

$$\varphi(n) = (p_1 - 1) \cdot (p_2 - 1), \tag{3.2b}$$

which from now on we replace, for factorization purposes, with:

$$\varphi_s(n) = (a - 1) \cdot (b - 1) \tag{3.3}$$

from which we obtain the following equation that defines the sum $\sum(a, b)$:

$$\sum(a, b) = n - \varphi_s(n) + 1, \tag{6.1b}$$

and with the sum $\sum(a, b)$ with which, denoting by $x_C$ and $y_C$ the respective values of $x$ and $y$ corresponding to the iteration with which we identify the first or unique pair of non-trivial factors of $n$, we can calculate the following value of $x = x_C$:

$$x_C = \frac{\sum(a, b)}{2} = \frac{a + b}{2} = \frac{n - \varphi_s(n) + 1}{2} \tag{6.30}$$

with which, by means of the well-known equation $y = \sqrt{x^2 - n}$ (or, resorting to the simplification described in the last part of Section 5 of this paper, i.e. resorting to (5.5), by means of the equation $y_r = \sqrt{y_f^2 + 4 \cdot (x_f + 1)}$, with $x_f = x_C - 2$), we obtain the integer value of $y = y_C$ (or the integer value of $y_r = y_C$) with which, together with the value of $x_C$, we identify the first or unique pair of non-trivial factors of $n$:

$$a = x_C + y_C \tag{6.31}$$

and

$$b = x_C - y_C. \tag{6.32}$$

Obviously, both $\varphi_s(n)$ and $\sum(a, b)$ are unknown a priori (with $a$ and $b$ being the factors of the first or unique pair of non-trivial factors of $n$); otherwise, the problem of identifying the first or unique pair of non-trivial factors of $n$ would be immediately solved.

Therefore, considering the fact that, according to the modification (for factorization purposes) of Euler's totient function $\varphi(n)$ of any $n$ that is a semiprime, both factors $a$ and $b$, with $a > b$, of any possible pair of non-trivial factors of $n$ and of the pair of trivial factors of $n$ are in any case always considered as prime even if they are not (therefore, both factors $a$ and $b$ of any pair of trivial factors of $n$, with $a$ which might not even be prime, are in any case always considered as prime, just as both factors $a$ and $b$ of any possible pair of non-trivial factors of $n$, which have prime factors in common and therefore are not coprime, are in any case always considered as prime), if $n$ is a non-semiprime composite and we decide to adopt DFM-1 up to the condition $x - y = 1$ (or, if we decide to resort to the simplification described in the last part of Section 5 of this paper by adopting (5.5), up to the condition $x_f + 2 - y_r = 1$), we will have the different values of $\varphi_s(n)$, associated with any pair of non-trivial factors $a$ and $b$ of $n$ and to the pair of trivial factors $a$ and $b$ of $n$, which will be always all different from the unique value of $\varphi(n)$ as determined due to the known assumptions of Euler's totient function $\varphi(n)$ [4] (highlighting that the value of $\varphi(n)$, associated with any $n$ that is a semiprime, is the product of the prime factors of $n$ which are decreased by 1). But, in such case, the values of $\varphi_s(n)$, associated with any pair of non-trivial factors of $n$ and to the pair of trivial factors of $n$, will be those effectively and validly recurring (with their validity which is absolutely confirmed by the results obtained with the application of DFM-1) and for which the equation $\sum(a, b) = n - \varphi_s(n) + 1$ is always respected, just as the equation $\sum(a, b) = n - \varphi(n) + 1$ is respected for the value of $\varphi(n)$ of any $n$ that is a semiprime. Indeed, if we consider the pairs of non-trivial factors and the pair of trivial factors of the examples represented in the previous Tables 1 and 2, we have, with the condition that any factor $a$ and $b$ of $n$ is in any case always considered as prime even if it is not, the different values of $\varphi_s(n)$, effectively and validly recurring during the application of DFM-1, as reported in Table 3 represented below:

Table 3. The different values of $\varphi_s(n)$ with $n = 70399$

| $n$ | pair of factors of $n$ | $a$ | $b$ | $\varphi_s(n) = (a-1) \cdot (b-1)$ | $\sum(a, b) = n - \varphi_s(n) + 1$ | $\varphi(n)$ |
|---|---|---|---|---|---|---|
| 70399 | 1st non-trivial | 623 | 113 | 69664 | 736 | $(113 - 1)$ $\cdot (89 - 1)$ $\cdot (7 - 1)$ $= 59136$ |
| | 2nd non-trivial | 791 | 89 | 69520 | 880 | |
| | 3rd non-trivial | 10057 | 7 | 60336 | 10064 | |
| | trivial | 70399 | 1 | 0 | 70400 | |

Obviously, if $n$ is a semiprime, we will have the value of $\varphi_s(n) = \varphi(n)$ that is associated with the unique pair of non-trivial factors $a$ and $b$ of $n$. And, obviously, if $n$ is a prime, we will reach the condition $x - y = 1$ (or $x_f + 2 - y_r = 1$), with $\varphi_s(n) = 0$ being $b - 1 = 0$, without having identified any pair of non-trivial factors $a$ and $b$ of $n$.

As is evident from the previous Table 3 and as already described in Lemma 7 and in the Theorem which have been stated in this section, as we identify the pairs of non-trivial factors $a$ and $b$ of $n$, the value of $a$ increases more and more and the value of $b$ decreases more and more, i.e. with $a$ and $b$ moving further and further away from each other, until the identification of the pair of trivial factors $a$ and $b$ of $n$, with the largest possible value of $a$ and the smallest possible value of $b$ (i.e. with $a$ and $b$ as far apart as possible), being $a = n$ and $b = 1$ (and with the value of $\varphi_s(n)$ decreasing more and more until the identification of the pair of trivial factors $a$ and $b$ of $n$ with the value of $\varphi_s(n) = 0$, and with the values of $\sum(a, b)$ increasing more and more until it reaches the value $\sum(a, b) = n + 1$, which is associated with the pair of trivial factors $a$ and $b$ of $n$). In this regard, I deem it appropriate to point out that, if we wanted to continue (even if without logical sense) incrementing by two units at a time the value of $x$ with which we identified the pair of trivial factors $a$ and $b$ of $n$, the value of $a$ will continue to increase more and more and the value of $b$ will continue to decrease more and more, i.e. with $a$ and $b$ that will continue to move further and further away from each other, and, obviously, with $a$ and $b$ that will always both be positive non-integers,(and with $b$ which will always be a positive number less than 1) but for which the following equation always holds:

$$\sum(a, b) = n - \varphi_{s1}(n) + (i - 1) \cdot 4 + 1, \qquad (6.33a)$$

i.e.:

$$\sum(a, b) = \sum(a_1, b_1) + (i - 1) \cdot 4, \qquad (6.33b)$$

and also:

$$\sum(a, b) = n - (a - 1) \cdot (b - 1) + 1, \qquad (6.33c)$$

with the value of $\sum(a, b)$ which will always be an even positive integer that will increase more and more, but with the product $(a - 1) \cdot (b - 1)$ which, in such infinite cases, will always be an even negative integer, a multiple of 4, with an absolute value that will increase more and more.

That said, and always having to consider that $n$ is a possible composite, we proceed in such a way that the computed initial value of $\varphi_s(n)$, which is $\varphi_{s1}(n)$ as indicated below:

$$\varphi_{s1}(n) = \left\lfloor \frac{\varphi_{s1_p}(n)}{4} \right\rfloor \cdot 4 = \left\lfloor \frac{n - 2 \cdot \lfloor \sqrt{n} \rfloor}{4} \right\rfloor \cdot 4, \qquad (6.13)$$

and which is therefore the largest possible value of $\varphi_s(n)$ (and with $\varphi_{s1_p}(n) = n - 2 \cdot \lfloor \sqrt{n} \rfloor$, which is the value of $\varphi_{s1}(n)$ provisional, which, as we know from Lemma 2, is an odd positive integer), comes to coincide with the value of $\varphi_s(n)$ associated with the first or unique pair of non-trivial factors of $n$, first checking whether the condition $\varphi_{s1}(n) = \varphi_s(n)$ occurs, and, in case of failure, decreasing the value of $\varphi_{s1}(n)$ by four units at a time. But, this action involves an unnecessary number of operations to be performed. So, to simplify the operations (reducing the number of operations to be carried out), we will proceed, after unsuccessfully checking for the condition $x_1 = x_C = \frac{n - \left\lfloor \frac{n - 2 \cdot \lfloor \sqrt{n} \rfloor}{4} \right\rfloor \cdot 4 + 1}{2}$, i.e. the condition $\varphi_{s1}(n) = \varphi_s(n)$, to increase the value of $x_1$, indicated below, always by two units at a time (starting, therefore, from the 2nd iteration with the value of $x_2 = x_1 + 2 = \frac{n - \left\lfloor \frac{n - 2 \cdot \lfloor \sqrt{n} \rfloor}{4} \right\rfloor \cdot 4 + 1}{2} + 2$) until reaching the condition $x = x_C$ using, for each iteration, the well-known equation $y = \sqrt{x^2 - n}$ (i.e. the equation $y^2 = x^2 - n$):

$$x_1 = \frac{n - \varphi_{s1}(n) + 1}{2}, \tag{6.14}$$

with the increase, always by two units at a time in the value of $x_1$ (or of $x_f = x_1$), which is evident from the simplification given in Section 6 with (6.18). Or (resorting to the simplification described in the last part of Section 5 of this paper and starting, therefore, from the 2nd iteration with the value of $x_f = x_1 = \frac{n - \left\lfloor \frac{n - 2 \cdot \lfloor \sqrt{n} \rfloor}{4} \right\rfloor \cdot 4 + 1}{2}$), we will proceed—after unsuccessfully checking for the condition $x_1 = x_C = \frac{n - \left\lfloor \frac{n - 2 \cdot \lfloor \sqrt{n} \rfloor}{4} \right\rfloor \cdot 4 + 1}{2}$, i.e. the condition $\varphi_{s1}(n) = \varphi_s(n)$—to increase the value of $x_f = x_1$ always by two units at a time until reaching the condition $x_f + 2 = x_C$ (where the value of $x_f$, in this case, is the value of $x$ corresponding to the iteration preceding the one in which we have identified the integer value of $y_r = y_C$, i.e. the integer value of $y_r$ with $x_f + 2 - y_r \neq 1$), using, for each iteration subsequent to the 1st iteration, the equation $y_r^2 = y_f^2 + 4 \cdot (x_f + 1)$ as given in Section 5 with (5.5).

Furthermore, the simplification, given in Section 6 with (6.18), facilitates the understanding of the circumstances regarding the adopted step-size factor and the hypothetical scenario, circumstances that we will present shortly.

It should be noted that, in the cases of $n = 7$, $n = 5$ and $n = 3$, the value of $\varphi_{s1}(n)$—as well as the value of $\varphi_s(n)$ associated with the pair of trivial factors $a$ and $b$ of any $n$ (being $a = n$ and $b = 1$)—is 0, and, therefore, we have the condition $\varphi_s(n) = \varphi_{s1}(n) = 0$.

Furthermore, we already know that, in the equation that defines the value of $x_1$, we have considered the integer part of $\frac{\varphi_{s1_p}(n)}{4}$, i.e. the floor of $\frac{\varphi_{s1_p}(n)}{4}$; and this since, obviously, the value of $\varphi_{s1}(n)$, i.e. $\left\lfloor \frac{\varphi_{s1_p}(n)}{4} \right\rfloor \cdot 4$, must be an integer and a multiple of 4 and, therefore, divisible by 4.

And we already know that, still in the equation that defines the value of $x_1$, we have considered the value of $\lfloor \sqrt{n} \rfloor$ which is the integer part of $\sqrt{n}$, i.e. the floor of $\sqrt{n}$. The value of $\lfloor \sqrt{n} \rfloor$ is the right value to consider to correctly calculate the value of $x_1$. Indeed, if we also adopt the decimal part of $\sqrt{n}$, in many cases, i.e. in infinite cases, we will never be able to identify the first or unique pair of non-trivial factors of $n$. These are the cases, which are infinite, among those for which, if, with the adoption of FFA-1, the first or unique pair of non-trivial factors of $n$ is identified at the 1st iteration, the value of $\varphi_{s1_p}(n)$, calculated without the decimal part of $\sqrt{n}$, will always turn out to be one unit greater than a number divisible by 4, i.e. one unit greater than the value of $\varphi_s(n)$ and such as to determine the condition $\varphi_{s1}(n) = \varphi_s(n)$; for this reason, if, in such cases, the value of $\sqrt{n}$ is given by the sum of the floor of $\sqrt{n}$ and a decimal number greater than 0.5, the value of $\varphi_{s1_p}(n)$, calculated with the decimal part of $\sqrt{n}$, will be lower than the value of $\varphi_{s1_p}(n)$, calculated without the decimal part of $\sqrt{n}$, of a value certainly greater than 1, thereby determining a value of $\varphi_{s1}(n)$ of four units lower than the aforementioned number divisible by 4, i.e. of four units lower than the value of $\varphi_s(n)$, and, consequently, a value of $x_1$ greater than one step size, i.e. greater than two units, compared to the value of $x_1$ with which we can identify the first or unique pair of non-trivial factors of n and, therefore, we will never be able to identify the first or unique pair of non-trivial factors of $n$.

And, obviously, if, with the adoption of FFA-1, the first or unique pair of non-trivial factors of $n$ is not identified at the 1st iteration and if the aforementioned two conditions occur (with the value of $\varphi_{s1_p}(n)$, calculated without the decimal part of $\sqrt{n}$, which is one unit greater than a number divisible by 4, and with the value of $\sqrt{n}$ given by the sum of the floor of $\sqrt{n}$ and a decimal number greater than 0.5), it happens that, starting from the value of $x_1$ calculated with the decimal part of $\sqrt{n}$, DFM-1 will be an improvement of one iteration compared to the adoption of the same DFM-1 starting from the value of $x_1$ calculated without the decimal part of $\sqrt{n}$.

While, if the value of $\varphi_{s1_p}(n)$, calculated without the decimal part of $\sqrt{n}$, will be three units greater than a number divisible by 4 (being $\varphi_{s1_p}(n)$, calculated without the decimal part of $\sqrt{n}$, always an odd integer), the value of $\varphi_{s1}(n)$, calculated with the decimal part of $\sqrt{n}$, will always be equal to the value of $\varphi_{s1}(n)$ calculated without the decimal part of $\sqrt{n}$, being the difference between the respective values of $\varphi_{s1_p}(n)$ always less than 2, and, consequently, the value of $x_1$, calculated with the part decimal part of $\sqrt{n}$, will always be equal to the value of $x_1$ calculated without the decimal part of $\sqrt{n}$. This circumstance always occurs in cases in which, with the adoption of FFA-1, the first or unique pair of non-trivial factors of $n$ is identified at the 2nd iteration and with the value of $\varphi_{s1_p}(n)$, calculated without the decimal part of $\sqrt{n}$, which is always three units greater than the value of $\varphi_s(n)$ and such to be determined, obviously, always and in any case the condition $\varphi_{s1}(n) = \varphi_s(n)$; for this reason, with DFM-1 and starting from the value of $x_1$ calculated with the decimal part of $\sqrt{n}$, an improvement of one iteration, compared to the adoption of the same DFM-1 starting from the value of $x_1$ calculated without the decimal part of $\sqrt{n}$, can only occur starting from the cases, which are infinite and in which the aforementioned two conditions occur (with the value of $\varphi_{s1_p}(n)$, calculated without the decimal part of $\sqrt{n}$, which is one unit greater than a number divisible by 4, and with the value of $\sqrt{n}$ given by the sum of the floor of $\sqrt{n}$ and a decimal number greater than 0.5), among those for which, with the adoption of FFA-1, the first or unique pair of non-trivial factors of $n$ is identified at the 3rd iteration and with the value of $\varphi_{s1_p}(n)$, calculated without the decimal part of $\sqrt{n}$, such as to always determine a value of $\varphi_{s1}(n)$ greater than four units compared to the value of $\varphi_s(n)$ and, consequently, a value of $x_1$ less than one step size, i.e. less than two units, compared to the value of $x_1$, calculated with the part decimal part of $\sqrt{n}$, with which we identify the first or unique pair of non-trivial factors of $n$.

And, again obviously, we can already easily know, in advance, if the aforementioned two conditions occur and if, with the adoption of FFA-1, the first or unique pair of non-trivial factors of $n$ cannot be identified at the 1st iteration (and already knowing a priori that, in the case in which, with the adoption of FFA-1, the first or unique pair of non-trivial factors of $n$ were to be identified at the 2nd iteration, the first of the aforementioned two conditions can never occur), so as to have to decide, possibly, to start from the value of $x_1$ calculated with the decimal part of $\sqrt{n}$ in order to reduce the number of iterations by one unit; but, we will have unnecessarily done extra operations to reduce the number of iterations by just one unit.

Therefore, the choice to start from the value of $x_1$ always calculated without the decimal part of $\sqrt{n}$ is quite justified.

That said, the number 4, which appears in the calculation of the value of $\varphi_{s1}(n)$—see (6.13)—both as the denominator of the fraction $\frac{\varphi_{s1_p}(n)}{4}$ and as the multiplier of $\left\lfloor\frac{\varphi_{s1_p}(n)}{4}\right\rfloor$, is the product between 2, as the permanent factor that we denote by $u$, and 2 as the step-size factor that we denote by $s$, both of which are factors that certainly form, among other factors, the value of $\varphi_s(n)$ associated with the first or unique pair of non-trivial factors of $n$; certainly since it is certain that the value of $\varphi_s(n)$, associated with the first or unique pair of non-trivial factors of $n$, is formed, among other factors, by the factor 2 with an exponent at least equal to 2 (see Lemma 3), being even both $(a - 1)$ and $(b - 1)$ which form the value of $\varphi_s(n)$ associated with the first or unique pair of non-trivial factors of $n$ (just as it is certain that the values of $\varphi_s(n)$, associated with any further pairs of non-trivial factors of $n$, are formed, among other factors, by the factor 2 with an exponent at least equal to 2). Thus, if we knew a priori the value of $\varphi_s(n)$ associated with the first or unique pair of non-trivial factors of $n$, we could also know a priori, as shown below, the number of iterations, that we denote by $i_{C_D}$, required by DFM-1 to reach the condition $x = x_C$ (identification of the first or unique pair of non-trivial factors of $n$):

$$i_{C_D} = \frac{\varphi_{s1}(n) - \varphi_s(n)}{u \cdot s} + 1 = \frac{\varphi_{s1}(n) - \varphi_s(n)}{2 \cdot 2} + 1 = \frac{\varphi_{s1}(n) - \varphi_s(n)}{4} + 1 \tag{6.34}$$

and with the *unit* that is *added* since we proceed, first of all, with the *initial check* for

$$x_1 = \frac{n - \left\lfloor\frac{n - 2 \cdot \lfloor\sqrt{n}\rfloor}{u \cdot s}\right\rfloor \cdot u \cdot s + 1}{2} = \frac{n - \left\lfloor\frac{n - 2 \cdot \lfloor\sqrt{n}\rfloor}{2 \cdot 2}\right\rfloor \cdot 2 \cdot 2 + 1}{2}. \tag{6.35}$$

## 7. Greater Efficiency of DFM-1 Compared to FFA-1

The step-size factor $s$, as defined in the last paragraph of the previous Section 6, represents the number by which, with the adoption of DFM-1, we reduce, possibly rounding up to the next integer, the number of iterations required by FFA-1 [2; 3, p. 1470] to identify the first or unique pair of non-trivial factors of $n$ (with the reduction which, obviously, certainly does not exist in cases in which the value of the number of iterations, required by FFA-1 to identify the first or unique pair of non-trivial factors of $n$, is equal to 1). Specifically, denoting by $i_{C_F}$ the number of iterations required by FFA-1—see (4.3)—to identify the first or unique pair of non-trivial factors of $n$, we observe that, with $s = 2$, the value of $i_{C_D}$—see (5.3) and (6.34)—assumes, in relation to the value of $i_{C_F}$, one of the following two forms (with the value of $i_{C_D}$ which includes, obviously, the iteration for $x_1 = \frac{n - \left\lfloor \frac{n - 2 \cdot \lfloor \sqrt{n} \rfloor}{u \cdot s} \right\rfloor \cdot u \cdot s + 1}{2}$):

$$i_{C_D} = \frac{i_{C_F}}{s}, \tag{7.1a}$$

where $\frac{i_{C_F}}{s} = \frac{i_{C_F}}{2}$ is clearly an integer, or

$$i_{C_D} = \left\lfloor \frac{i_{C_F}}{s} \right\rfloor + 1 = \left\lceil \frac{i_{C_F}}{s} \right\rceil, \tag{7.1b}$$

where $\frac{i_{C_F}}{s} = \frac{i_{C_F}}{2}$ is not an integer.

In short, it is true that the asymptotic complexity of DFM-1 does not change with respect to that of FFA-1 (an asymptotic complexity that is $O(n)$ in the worst case), but it is also true that, to identify the first or unique pair of non-trivial factors of $n$, the computational cost $T(n)$ of DFM-1 is reduced by half or by approximately half (with the reduction which, obviously, certainly does not exist in cases in which the value of the number of iterations, required by FFA-1 to identify the first or unique pair of non-trivial factors of $n$, is equal to 1) compared to that of FFA-1, going from

$$T(n) = k \cdot M(\log n), \tag{7.2}$$

which is the computational cost $T(n)$ of FFA-1, to

$$T(n) = \frac{k}{2} \cdot M(\log n) \tag{7.3a}$$

or to

$$T(n) = \left\lceil \frac{k}{2} \right\rceil \cdot M(\log n), \tag{7.3b}$$

which are the computational costs $T(n)$ of DFM-1 (with $s = 2$), where $k = i_{C_F}$, $\log n$ is the number of bits required to represent the integer $n$ being factored, and $M(\log n)$ represents the computational cost of multiplying two numbers of size $\log n$. As is evident, these indications on the reduction of the computational cost consider the number of iterations to be performed, having $k = i_{C_F}$, and not the number of operations to be performed. This is due to the fact that, if we consider the number of operations to be performed and we resort (both with FFA-1 and with DFM-1) to the simplification that concerns the possible elimination of the square root operation by observing the last digit of $y^2 = x^2 - n$ (as well as the last digit of $y_r^2$ and as well as the last digit of $y_1^2$, where $y_1^2$ is the first value of $y_f^2$ possibly to be considered), the above reductions in computational cost, as the distance between the factors $a$ and $b$ of the first or unique pair of non-trivial factors of $n$ increases, decrease more and more, although they always remain relevant, since, with FFA-1, the cases in which the square root operation is eliminated increase more and more as we perform the iterations that don't occur with DFM-1. Therefore, in order to simplify the indications on the reduction of the computational cost of DFM-1 (with $s = 2$) compared to that of FFA-1, we will always and only consider the number of iterations to be performed and not the number of operations to be performed.

That said, and returning to the form that the value of $i_{C_D}$ assumes in relation to the value of $i_{C_F}$, if, for example, FFA-1 requires 56 iterations to identify the unique pair of non-trivial factors of the semiprime 8612553881 = 96059 · 89659, DFM-1 (with $s = 2$) requires only 28 iterations (including the iteration for $x_1 = \frac{n - \left\lfloor \frac{n - 2 \cdot \lfloor \sqrt{n} \rfloor}{u \cdot s} \right\rfloor \cdot u \cdot s + 1}{2}$) for the same identification and, therefore, we have the value of $i_{C_D}$ that is of the form $\frac{i_{C_F}}{s}$.

And, to give another example, if FFA-1 requires 1515 iterations to identify the unique pair of non-trivial factors of the semiprime 5357811983 = 89681 · 59743, DFM-1 (with $s = 2$) requires only 758 iterations (including the iteration for $x_1 = \frac{n - \left\lfloor \frac{n - 2 \cdot \lfloor \sqrt{n} \rfloor}{u \cdot s} \right\rfloor \cdot u \cdot s + 1}{2}$) for the same identification and, therefore, we have the value of $i_{C_D}$ that is of the form $\left\lfloor \frac{i_{C_F}}{s} \right\rfloor + 1 = \left\lceil \frac{i_{C_F}}{s} \right\rceil$.

Obviously, the aforementioned indications on the form that the value of $i_{C_D}$ assumes in relation to the value of $i_{C_F}$, and all the aspects indicated above concerning the computational cost $T(n)$ of DFM-1, also hold both for identifying any pair of non-trivial factors $a$ and $b$ of $n$ and for identifying the pair of trivial factors $a$ and $b$ of $n$, and they also hold to reach the limit number of iterations, to be specified later, which, if the conditions occur, may also be reached at the 1st iteration, i.e. with the value of $x_1 = \frac{n - \left\lfloor \frac{n - 2 \cdot \lfloor \sqrt{n} \rfloor}{u \cdot s} \right\rfloor \cdot u \cdot s + 1}{2} = \frac{n - \left\lfloor \frac{n - 2 \cdot \lfloor \sqrt{n} \rfloor}{2 \cdot 2} \right\rfloor \cdot 2 \cdot 2 + 1}{2}$.

Therefore, for the purposes of identifying the first or unique pair of non-trivial factors of $n$, both the permanent factor $u$ and the step-size factor $s$ must necessarily be factors that form, among other factors, the value of $\varphi_s(n)$ associated with the first or unique pair of non-trivial factors of $n$; and, since, as it has already been highlighted, this value of $\varphi_s(n)$ is certainly formed, among other factors, by the factor 2 with an exponent at least equal to 2 (obviously, as we already know, this is also true for the value of $\varphi_s(n)$ associated with any pair of non-trivial factors of $n$), with DFM-1 we can certainly proceed to halve, possibly rounding up to the next integer, the number of iterations required by FFA-1 to identify the first or unique pair of non-trivial factors of $n$ (highlighting that, with $i_{C_F} = 1$ and $s = 2$, obviously the value of $i_{C_D}$ is of the form $\left\lfloor \frac{i_{C_F}}{s} \right\rfloor + 1 = \left\lceil \frac{i_{C_F}}{s} \right\rceil$ being $i_{C_D} = \left\lfloor \frac{1}{2} \right\rfloor + 1 = \left\lceil \frac{1}{2} \right\rceil = 1$).

## 8. Hypothetical Scenario

At this point, it is necessary to highlight that, hypothesizing, reasonably (for reasons of probability), that the value of $\varphi_s(n)$, associated with the first or unique pair of non-trivial factors of $n$, can be formed, among other possible odd prime factors, by the factor 2 with an exponent greater than 2, with DFM-1 it could be possible (they are attempts), for example in the presence of a hypothetical and probable factor $2^3$ or $2^4$ or $2^5$, to reduce the number of iterations, required by FFA-1 to identify the first or unique pair of non-trivial factors of $n$, to 1/4 with $s = 2^2$ or to 1/8 with $s = 2^3$ or to 1/16 with $s = 2^4$ (always possibly rounding up to the next integer); and, at the same time, hypothesizing, reasonably (always for reasons of probability), that also the factor 3 occurs in the value of $\varphi_s(n)$ associated with the first or unique pair of non-trivial factors of $n$, with DFM-1 we could try to reduce the number of iterations, required by FFA-1 to identify the first or unique pair of non-trivial factors of $n$, to 1/3 with $s = 3$ or to 1/6 with $s = 2 \cdot 3$ or to 1/12 with $s = 2^2 \cdot 3$ or to 1/24 with $s = 2^3 \cdot 3$ or to 1/48 with $s = 2^4 \cdot 3$ (and always possibly rounding up to the next integer).

Therefore, in a hypothetical line, we can also proceed with a step-size factor $s$, which contributes to forming the value of $\varphi_s(n)$ associated with the first or unique pair of non-trivial factors of $n$, which is composite, just as we can also proceed, always in a hypothetical line, with a step-size factor $s$, which contributes to forming the value of $\varphi_s(n)$ associated with the first or unique pair of non-trivial factors of $n$, which is only a prime other than 2 or which is an odd composite. But, we will never adopt a step-size factor $s$ that is odd since we already know that the step-size factor $s$ is certainly always formed, among other possible odd prime factors, by the factor 2 with an exponent at least equal to 1. Obviously, in all cases, the permanent factor $u = 2$ remains unchanged, so that, in the calculation of the value of $\varphi_{s1}(n)$ [see (6.13)], the aforementioned number 4 (as specified in the last part of Section 6), both as the denominator of the fraction $\frac{\varphi_{s1_p}(n)}{4}$ and as the multiplier of $\left\lfloor \frac{\varphi_{s1_p}(n)}{4} \right\rfloor$, is replaced, respectively in the aforementioned hypotheses, by the numbers 8, 16 and 32 and by the numbers 6, 12, 24, 48 and 96. To give an example, if, considering the aforementioned hypotheses, we wanted to proceed with $s = 2^4 \cdot 3 = 48$, the value of $\varphi_{s1}(n)$ would be calculated as follows:

$$\varphi_{s1}(n) = \left\lfloor \frac{\varphi_{s1_p}(n)}{u \cdot s} \right\rfloor \cdot u \cdot s = \left\lfloor \frac{\varphi_{s1_p}(n)}{2 \cdot 2^4 \cdot 3} \right\rfloor \cdot 2 \cdot 2^4 \cdot 3 = \left\lfloor \frac{n - 2 \cdot \lfloor \sqrt{n} \rfloor}{96} \right\rfloor \cdot 96.$$

As is evident, the expression for calculating the value of $\varphi_{s1}(n)$ provisional, which is $\varphi_{s1_p}(n)$, always remains unchanged for any step size adopted.

And, obviously, also for the aforementioned hypotheses, should they correspond to truly recurring circumstances, the value of $i_{C_D}$, required by DFM-1 to identify the first or unique pair of non-trivial factors of $n$, will always assume, in relation to the value of $i_{C_F}$ and in relation to the step-size factor $s$ adopted, one of the following two forms as indicated in the previous Section 7 with $s = 2$ (and with the value of $i_{C_D}$ which includes, obviously, the iteration for $x_1 = \frac{n - \left\lfloor \frac{n - 2 \cdot \lfloor \sqrt{n} \rfloor}{u \cdot s} \right\rfloor \cdot u \cdot s + 1}{2}$):

$$i_{C_D} = \frac{i_{C_F}}{s}, \tag{7.1$a$}$$

with $\frac{i_{C_F}}{s}$ which is obviously an integer, or

$$i_{C_D} = \left\lfloor \frac{i_{C_F}}{s} \right\rfloor + 1 = \left\lceil \frac{i_{C_F}}{s} \right\rceil, \tag{7.1$b$}$$

with $\frac{i_{C_F}}{s}$ which is a non-integer, and, therefore, as regards the example represented above with $s = 2^4 \cdot 3 = 48$, and provided that the hypothesis in question corresponds to a truly recurring circumstance, the value of $i_{C_D}$, required by DFM-1 to identify the first or unique pair of non-trivial factors of $n$, will assume, in relation to the value of $i_{C_F}$, one of the following two forms (with the value of $i_{C_D}$ which includes, obviously, the iteration for $x_1 = \frac{n - \left\lfloor \frac{n - 2 \cdot \lfloor\sqrt{n}\rfloor}{2 \cdot 48} \right\rfloor \cdot 2 \cdot 48 + 1}{2}$):

$$i_{C_D} = \frac{i_{C_F}}{s} = \frac{i_{C_F}}{48},$$

with $\frac{i_{C_F}}{s} = \frac{i_{C_F}}{48}$ which is obviously an integer, or

$$i_{C_D} = \left\lfloor \frac{i_{C_F}}{s} \right\rfloor + 1 = \left\lfloor \frac{i_{C_F}}{48} \right\rfloor + 1 = \left\lceil \frac{i_{C_F}}{s} \right\rceil = \left\lceil \frac{i_{C_F}}{48} \right\rceil,$$

with $\frac{i_{C_F}}{s} = \frac{i_{C_F}}{48}$ which is a non-integer.

Just as, also for the aforementioned hypotheses, the above indications on the form that the value of $i_{C_D}$ assumes in relation to the value of $i_{C_F}$ also hold, but always in relation to the step-size factor $s$ adopted, both for identifying any pair of non-trivial factors $a$ and $b$ of $n$ and for identifying the pair of trivial factors $a$ and $b$ of $n$, and they also hold to reach the limit number of iterations, to be specified later, which, if the conditions occur, may also be reached at the 1st iteration, i.e. with the value of $x_1 = \frac{n - \left\lfloor \frac{n - 2 \cdot \lfloor\sqrt{n}\rfloor}{u \cdot s} \right\rfloor \cdot u \cdot s + 1}{2}$. But, keeping always in mind that, by adopting a step-size factor $s$ greater than 2, we will never be certain of the primality of $n$, even though we identify the pair of trivial factors of $n$ or we reach the limit number of iterations without having identified any pair of non-trivial factors of $n$.

However, it is worth highlighting the circumstance for which, if we adopted a step-size factor $s$, which contributes to forming the value of $\varphi_s(n)$ associated with the first or unique pair of non-trivial factors of $n$, with a value equal to or greater than the value of the number of iterations required by FFA-1 to identify the first or unique pair of non-trivial factors of $n$ (but, obviously, we cannot know this a priori), we would identify, with the adoption of DFM-1, the first or unique pair of non-trivial factors of $n$ always with a single iteration, i.e. always at the 1st iteration and, therefore, we will always have the condition $x_1 = x_C$. Obviously, this is true for identifying any pair of non-trivial factors $a$ and $b$ of $n$, i.e. with the adoption of DFM-1 we would identify a specific pair of non-trivial factors of $n$ always with a single iteration, i.e. always at the 1st iteration with $y_1$ which is an integer and with $x_1 - y_1 \neq 1$, by adopting a step-size factor $s$, which contributes to forming the value of $\varphi_s(n)$ associated with the specific pair of non-trivial factors of $n$, with a value equal to or greater than the value of the number of iterations required by FFA-1 to identify the specific pair of non-trivial factors of $n$ (obviously, there cannot exist a step-size factor $s$ which, as such, contributes to forming the value of $\varphi_s(n)$ associated with two different pairs of non-trivial factors of $n$ and that has, at the same time, both a value equal to or greater than the value of the number of iterations required by FFA-1 to identify one of the two different pairs of non-trivial factors of $n$, and a value equal to or greater than the value of the number of iterations required by FFA-1 to identify the other of the two different pairs of non-trivial factors of $n$). While and for the reasons that will be more evident later, to identify, with the adoption of DFM-1, the pair of trivial factors of $n$ (a pair which, however, we already know a priori, since the corresponding values of $x = \frac{n+1}{2}$ and $y = \frac{n-1}{2}$ are already known a priori) always with a single iteration, i.e. always at the 1st iteration with $y_1 = \frac{n-1}{2}$ which is obviously an integer, we can adopt any step-size factor $s$ with a value equal to or greater than the value of the number of iterations required by FFA-1 to identify the pair of trivial factors of $n$.

Now, still remaining in the hypothetical scenario and before proceeding to represent the phases of the new factorization algorithm (DFM-1) in question (with $s = 2$), it is also necessary to highlight that, with $n$ for which more than one pair of non-trivial factors occurs (but, we cannot know a priori), it may happen, always without prejudice to the permanent factor $u = 2$, that we identify the pair of non-trivial factors of $n$ which is not the first pair of non-trivial factors of $n$. This happens when, always by hypothesis, we adopt a step-size factor $s$ that contributes to forming the value of $\varphi_s(n)$ associated with the identified pair of non-trivial factors of $n$, but which does not contribute to forming the value of $\varphi_s(n)$ associated with the first pair of non-trivial factors of $n$.

But, as evident is, in all the hypothetical cases represented above we therefore enter the field of hypotheses, although the hypotheses may refer, with great probability, to truly recurring circumstances (if we consider the fact that small prime factors, such as the primes 2 and 3, can form the value of $\varphi_s(n)$, associated with any pair of non-trivial factors of $n$, with a higher probability).

In Table 4 represented below and always considering the example with $n = 70399$, we compare, with regard to the identification of the first pair of non-trivial factors of $n$, FFA-1, DFM-1 with $s = 2$, and DFM-1 with hypothetical step-size factor $s = 2^3$, indicating, in addition to the values of $\varphi_{s1}(n)$ and $x_1$, the values of $i_{C_F}$ and $i_{C_D}$ associated with the first pair of non-trivial factors of $n$, and with $\varphi_s(n) = (a-1) \cdot (b-1) = (623-1) \cdot (113-1) = 69664 = 2^5 \cdot 7 \cdot 311$, which is associated with the first pair of non-trivial factors of $n$, which, obviously , we do not know a priori:

Table 4. The comparison of three algorithms for identifying the first pair of non-trivial factors of $n = 70399$

| algorithm | $\varphi_{s1}(n)$ | $x_1$ | $i_{C_F}$ | $i_{C_D}$ |
|---|---|---|---|---|
| FFA-1 (with $s = 1$) | 69868 | 266 | 103 | $103 = \frac{i_{C_F}}{1}$ |
| DFM-1 with $s = 2$ | 69868 | 266 | | $52 = \left\lfloor \frac{i_{C_F}}{2} \right\rfloor + 1 = \left\lceil \frac{i_{C_F}}{2} \right\rceil$ |
| DFM-1 with hypothetical step size $s = 2^3$ | 69856 | 272 | | $13 = \left\lfloor \frac{i_{C_F}}{2^3} \right\rfloor + 1 = \left\lceil \frac{i_{C_F}}{2^3} \right\rceil$ |

As we can see, in Table 4 we reported, in the line concerning FFA-1, the values of $\varphi_{s1}(n)$ and $i_{C_D}$; and this is since, with DFM-1, we can obviously also proceed with a step-size factor $s = 1$, with the respective value of $x_1 = \frac{n - \left\lfloor \frac{n - 2 \cdot \lfloor \sqrt{n} \rfloor}{u \cdot s} \right\rfloor \cdot u \cdot s + 1}{2} = \frac{n - \left\lfloor \frac{n - 2 \cdot \lfloor \sqrt{n} \rfloor}{2 \cdot 1} \right\rfloor \cdot 2 \cdot 1 + 1}{2}$ that, contrary to what happens when adopting a step-size factor $s > 1$, is always equal to $\lceil \sqrt{n} \rceil$.

That said, we conclude this section by highlighting that, if we consider the first example represented in Section 7 of this paper, with $n = 8612553881$ and with the value of $\varphi_s(n) = 8612368164 = 2^2 \cdot 3^2 \cdot 17 \cdot 293 \cdot 48029$ which is associated with the unique pair of non-trivial factors of $n$, by adopting DFM-1 with a hypothetical step-size factor $s = 2 \cdot 3 = 6$ (which contributes to forming the value of $\varphi_s(n)$ associated with the unique pair of non-trivial factors of $n$, but we cannot know a priori, although the probability that this step-size factor could actually occur is very great), only 10 iterations (including the iteration for $x_1 = \frac{n - \left\lfloor \frac{n - 2 \cdot \lfloor \sqrt{n} \rfloor}{u \cdot s} \right\rfloor \cdot u \cdot s + 1}{2} = \frac{n - \left\lfloor \frac{n - 2 \cdot \lfloor \sqrt{n} \rfloor}{2 \cdot 6} \right\rfloor \cdot 2 \cdot 6 + 1}{2}$) are required for the identification of the unique pair of non-trivial factors of $n$, instead of 56 iterations required by FFA-1 for the same identification, and, therefore, we have the value of $i_{C_D}$ that is of the form $\left\lfloor \frac{i_{C_F}}{s} \right\rfloor + 1 = \left\lfloor \frac{i_{C_F}}{2 \cdot 3} \right\rfloor + 1 = \left\lceil \frac{i_{C_F}}{s} \right\rceil = \left\lceil \frac{i_{C_F}}{2 \cdot 3} \right\rceil$. Instead, if we consider the second example always represented in Section 7 of this paper, with $n = 5357811983$ and with the value of $\varphi_s(n) = 5357662560 = 2^5 \cdot 3^2 \cdot 5 \cdot 19 \cdot 59 \cdot 3319$ which is associated with the unique pair of non-trivial factors of $n$, by adopting DFM-1 with a hypothetical step-size factor $s = 2^2 \cdot 3 = 12$ (which contributes to forming the value of $\varphi_s(n)$ associated with the unique pair of non-trivial factors of $n$, but we cannot know a priori, although the probability that this step-size factor could actually occur is very great), only 127 iterations (including the iteration for $x_1 = \frac{n - \left\lfloor \frac{n - 2 \cdot \lfloor \sqrt{n} \rfloor}{u \cdot s} \right\rfloor \cdot u \cdot s + 1}{2} = \frac{n - \left\lfloor \frac{n - 2 \cdot \lfloor \sqrt{n} \rfloor}{2 \cdot 12} \right\rfloor \cdot 2 \cdot 12 + 1}{2}$) are required for the identification of the unique pair of non-trivial factors of $n$, instead of 1515 iterations required by FFA-1 for the same identification, and, therefore, we have the value of $i_{C_D}$ that is of the form $\left\lfloor \frac{i_{C_F}}{s} \right\rfloor + 1 = \left\lfloor \frac{i_{C_F}}{2^2 \cdot 3} \right\rfloor + 1 = \left\lceil \frac{i_{C_F}}{s} \right\rceil = \left\lceil \frac{i_{C_F}}{2^2 \cdot 3} \right\rceil$.

## 9. Phases of DFM-1 (With s = 2)

Therefore, it is certain only the halving (possibly rounding up to the next integer) of the number of iterations required by FFA-1 to identify any pair of non-trivial factors $a$ and $b$ of $n$ and to identify the pair of trivial factors $a$ and $b$ of $n$, and also to reach the limit number of iterations, to be specified later, which, if the conditions occur, may also be reached at the 1st iteration, i.e. with the value of $x_1 = \frac{n - \left\lfloor \frac{n - 2 \cdot \lfloor \sqrt{n} \rfloor}{2 \cdot 2} \right\rfloor \cdot 2 \cdot 2 + 1}{2}$ (returning to highlight that, with $i_{C_F} = 1$ and $s = 2$, obviously the value of $i_{C_D}$ is of the form $\left\lfloor \frac{i_{C_F}}{s} \right\rfloor + 1 = \left\lceil \frac{i_{C_F}}{s} \right\rceil$ being $i_{C_D} = \left\lfloor \frac{1}{2} \right\rfloor + 1 = \left\lceil \frac{1}{2} \right\rceil = 1$). And therefore, it is certain the adoption of DFM-1 with $s = 2$.

Obviously, the adoption of DFM-1 with $s = 1$ is also certain, but without obtaining any advantage compared to FFA-1, highlighting that, with $s = 1$, obviously the value of $i_{C_D}$ is always of the form $\frac{i_{C_F}}{s}$. In this regard, It is appropriate to highlight that (4.3), i.e. the equation $i = \frac{a+b}{2} - \lfloor \sqrt{n} \rfloor$, can therefore also be written as follows, with $s = 1$:

$$i_{C_F} = \frac{\frac{a+b}{2} - \lceil \sqrt{n} \rceil}{s} + 1. \tag{9.1}$$

That said, below we present the phases of the new factorization algorithm (DFM-1) in question (with $s = 2$), which will be preceded by the calculation of the following initial value of $x = x_1$ [see (5.1), (6.7), (6.13), (6.14), and (6.35)]:

$$x_1 = \frac{n - \varphi_{s1}(n) + 1}{2}$$

$$= \frac{n - \left\lfloor \frac{\varphi_{s1_p}(n)}{u \cdot s} \right\rfloor \cdot u \cdot s + 1}{2}$$

$$= \frac{n - \left\lfloor \frac{n - 2 \cdot \lfloor \sqrt{n} \rfloor}{2 \cdot 2} \right\rfloor \cdot 2 \cdot 2 + 1}{2}. \tag{9.2}$$

Therefore, the phases of DFM-1 are as follows:

- **1st phase** - Check of the values of $y = y_1$ for $x = x_1$ [see (2.7), (6.15), and (6.17)]:

$$y_1 = \sqrt{x_1{}^2 - n}; \tag{9.3}$$

if the value of $y = y_1$ is an integer (with $x_1 - y_1 \neq 1$), we will proceed with the 2nd phase:

**2nd phase** - Calculation of the values of $a$ and $b$ [see (2.9), (2.10), (6.28), (6.29), (6.31), and (6.32)]:

$$a = x_1 + y_1 = x_C + y_C \tag{9.4}$$

and

$$b = x_1 - y_1 = x_C - y_C, \tag{9.5}$$

thus identifying, at the 1st iteration, the first or unique pair of non-trivial factors $a$ and $b$ of $n$.

While, if the value of $y = y_1$ is a non-integer, i.e. if the value of $y^2 = y_1{}^2$ is a non-perfect square, we will repeat the 1st phase, increasing the value of $x$ always by two units at a time (starting, therefore, from the 2nd iteration with the value of $x_2 = x_1 + 2 = \frac{n - \left\lfloor \frac{n - 2 \cdot \lfloor \sqrt{n} \rfloor}{4} \right\rfloor \cdot 4 + 1}{2} + 2$) in the equation $y = \sqrt{x^2 - n}$ (i.e. in the equation $y^2 = x^2 - n$), until we have identified the integer value of $y$ with $x - y \neq 1$, which is $y_C$, with which, together with the value of $x_C$, we identify, with the 2nd phase and as indicated below, the first or unique pair of non-trivial factors $a$ and $b$ of $n$:

$$a = x + y = x_C + y_C \tag{9.6}$$

and

$$b = x - y = x_C - y_C, \tag{9.7}$$

or until we are certain that an integer value of $y$ with $x - y \neq 1$can never occur such as to be able to identify the first or unique pair of non-trivial factors $a$ and $b$ of $n$ and, therefore, until we are certain of the primality of $n$.

Or—wishing to resort to the simplification described in the last part of Section 5 of this paper, and still in the case where the value of $y = y_1$ is a non-integer, i.e. still in the case where the value of $y^2 = y_1{}^2$ is a non-perfect square—we will adopt the equation $y_r{}^2 = y_f^2 + 4 \cdot (x_f + 1)$, as given in Section 5 with (5.5), for each iteration subsequent to the 1st iteration (starting, therefore, from the 2nd iteration with the value of $x_f = x_1 = \frac{n - \left\lfloor \frac{n - 2 \cdot \lfloor \sqrt{n} \rfloor}{4} \right\rfloor \cdot 4 + 1}{2}$ and, thereafter, increasing the value of $x_f$ always by two units at a time), until we have identified the integer value of of $y_r$ with $x_f + 2 - y_r \neq 1$ (where $x_f$ is the value of $x$ corresponding to the iteration preceding the one in which we have identified the integer value of $y_r$ with $x_f + 2 - y_r \neq 1$), which is $y_C$, with which, together with the value of $x_C$, we identify, with the 2nd phase and as indicated below, the first or unique pair of non-trivial factors $a$ and $b$ of $n$:

$$a = (x_f + 2) + y_r = x_C + y_C \tag{9.8}$$

and

$$b = (x_f + 2) - y_r = x_C - y_C, \tag{9.9}$$

or until we are certain that an integer value of $y_r$ with $x_f + 2 - y_r \neq 1$—with the values of $x_f$ and $y_r$ that appear in the same equation $y_r{}^2 = y_f^2 + 4 \cdot (x_f + 1)$—can never occur such as to be able to identify the first or unique pair of non-trivial factors $a$ and $b$ of $n$ and, therefore, until we are certain of the primality of $n$.

Obviously, DFM-1 (with $s = 2$) turns out to be of great importance, compared to FFA-1, in the presence of composites for which the identification of the first or unique pair of non-trivial factors of $n$, resorting to the adoption of FFA-1, requires a very large number of iterations. Indeed, considering, for example, very large semiprimes for which the identification of the unique pair of non-trivial factors of $n$, resorting to the adoption of FFA-1, requires $10^{21}$ iterations, the reduction of $5 \cdot 10^{20}$ iterations is a very significant result.

## 10. Limit number of iterations

At this point, and before proceeding with one final, detailed example of the application of DFM-1 (with $s = 2$), it is necessary to make a consideration on the identification, always with DFM-1 (with $s = 2$), of the pair of trivial factors of $n$ for the purposes of identifying the first or unique pair of non-trivial factors of $n$. It being understood that we already know a priori the pair of trivial factors $a$ and $b$ of any $n$, with $a = n$ and $b = 1$ and with the following corresponding values of $x$ and $y$, denoted by $x_P$ and $y_P$, respectively:

$$x_P = \frac{n + 1}{2} \tag{10.1}$$

and

$$y_P = \frac{n - 1}{2}, \tag{10.2}$$

it is obvious that it will always be necessary to carry out all the steps, starting from the value of $x = x_1 = \frac{n - \left\lfloor \frac{n - 2 \cdot \lfloor \sqrt{n} \rfloor}{4} \right\rfloor \cdot 4 + 1}{2}$ and also resorting to the simplifications described in the last part of Section 5 of this paper, and until we are certain that an integer value of $y$ with $x - y \neq 1$—or an integer value of $y_r$ with $x_f + 2 - y_r \neq 1$, with the values of $x_f$ and $y_r$ that appear in the same equation $y_r^2 = y_f^2 + 4 \cdot (x_f + 1)$—can never occur such as to be able to identify the first or unique pair of non-trivial factors of $n$ and, therefore, until we are certain of the primality of $n$.

But (wanting to open a parenthesis), it is also necessary to highlight that, unlike the identification of the first or unique pair of non-trivial factors of $n$ for which DFM-1 (even with a hypothetical step-size factor $s$ and, therefore, also with a step-size factor $s$ greater than 2) requires that the step-size factor $s$ adopted (like the permanent factor $u$) contributes to forming the value of $\varphi_s(n)$ associated with the first or unique pair of non-trivial factors of $n$ (this necessity, obviously and as we already know, is fundamental for the identification of each pair of non-trivial factors of $n$, with the respective value of $\varphi_s(n)$ which must necessarily be formed, among other possible factors, by the step-size factor $s$ adopted), the identification of the pair of trivial factors of $n$ occurs for any step-size factor $s$ adopted, therefore both with a step-size factor $s$ that contributes to forming the value of $\varphi_s(n)$ associated with the first or unique pair of non-trivial factors of $n$ and with a step-size factor $s$ that does not contribute to forming the value of $\varphi_s(n)$ associated with the first or unique pair of non-trivial factors of $n$; and this is due to the fact (which we already know from Lemma 4 stated in Section 6 of this paper) that, being the pair of trivial factors of $n$, the value of $\varphi_s(n)$ is 0, so much so that the value of $x = x_P$, to identify the pair of trivial factors of $n$, is exactly [see (6.30)]:

$$x_P = \frac{n - \varphi(n) + 1}{2} = \frac{n - 0 + 1}{2} = \frac{n + 1}{2}. \tag{10.3}$$

Therefore, if the step-size factor $s$ adopted does not contribute to forming the value of $\varphi_s(n)$ associated with the first or unique pair of non-trivial factors of $n$, with DFM-1 (even with a hypothetical step-size factor $s$ and, therefore, also with a step-size factor $s$ greater than 2) we will however and always arrive at the iteration with which we identify the pair of trivial factors of $n$, but without having identified the first or unique pair of non-trivial factors of $n$ (but, possibly, if further pairs of non-trivial factors of $n$ occur, such pairs can be identified if the respective values of $\varphi_s(n)$ are formed, among other factors, by the step-size factor $s$ adopted which does not form, among other factors, the value of $\varphi_s(n)$ associated with the first or unique pair of non-trivial factors of $n$).

But, obviously, with DFM-1 (with $s = 2$), it is not always necessary to reach the condition $x - y = 1$ (or the condition $x_f + 2 - y_r = 1$), with the following respective number of iterations (including the iteration for $x_1$) which we denote by $i_P$ [see (5.4) and (6.34)]:

$$\begin{aligned} i_P &= \frac{\varphi_{s1}(n) - \varphi_s(n)}{u \cdot s} + 1 \\ &= \frac{\varphi_{s1}(n) - \varphi_s(n)}{2 \cdot 2} + 1 \\ &= \frac{\varphi_{s1}(n) - 0}{2 \cdot 2} + 1 \\ &= \frac{\varphi_{s1}(n)}{2 \cdot 2} + 1; \end{aligned} \tag{10.4a}$$

or

$$i_P = \frac{\frac{n+1}{2} - x_1}{s} + 1 = \frac{\frac{n+1}{2} - x_1}{2} + 1 \qquad (10.4b)$$

or

$$i_P = \left\lfloor \frac{\varphi_{s1_p}(n)}{u \cdot s} \right\rfloor + 1 = \left\lfloor \frac{\varphi_{s1_p}(n)}{2 \cdot 2} \right\rfloor + 1 = \left\lceil \frac{\varphi_{s1_p}(n)}{2 \cdot 2} \right\rceil, \qquad (10.4c)$$

to have the certainty that an integer value of $y$ with $x - y \neq 1$ can never occur such as to be able to identify the first or unique pair of non-trivial factors of $n$ and, therefore, to be certain of the primality of $n$ (and, obviously and as already highlighted in the Introduction of the book *I moltiplicatori complessi per la possibile scomposizione dei Numeri RSA* [6, pp. 24-25], and as we can already deduce from the oft-cited reduction, with respect to FFA-1, also of the limit number of iterations, this “non-necessity”, obviously, also holds for the adoption of FFA-1, with the respective value of $i_P = \frac{n+1}{2} - \lfloor\sqrt{n}\rfloor$). Indeed, with DFM-1 (with $s = 2$), to have such certainty it will be sufficient to reach the limit number of iterations (which can be calculated a priori, as for the value of $i_P$), which we denote by $i_L$ and whose value is less than the respective value of $i_P$, without having identified an integer value of $y$, with $x - y \neq 1$, such as to be able to identify the first or unique pair of non-trivial factors of $n$.

As an absolute definition, the limit number of iterations $i_L$, which is less than the number of iterations required by DFM-1 to reach the condition $x - y = 1$—or the condition $x_f + 2 - y_r = 1$, with the values of $x_f$ and $y_r$ that appear in the same equation $y_r^2 = y_f^2 + 4 \cdot (x_f + 1)$—which characterizes the pair of trivial factors of $n$, represents the maximum number of iterations beyond which it is no longer possible for pairs of non-trivial factors of $n$ to occur.

Denoting by $b_L$ the smallest possible value to consider of b non-trivial of $n$ (thus, $b_L$ is a prime), with the limit number of iterations $i_L$ we identify the pair of values $x$ and $y$ that determines, on the basis of the value of $b_L$ adopted, the value of $x - y = b_L$, if the value of $b_L$ is equal to the smallest value (which may also be the only value if $n$ is a semiprime) of $b$ non-trivial of $n$, or the smallest value of $x - y$ greater than $b_L$ (with $y$ which is a not integer) if the value of $b_L$ is not equal to the smallest value (which may also be the only value if $n$ is a semiprime) of $b$ non-trivial of $n$ (but, it may also not occur any non-trivial value of $b$ of $n$ if $n$ is a prime). In particular, the value of $b_L$ to be adopted can never exceed the maximum value, which we denote by $b_{L_{max}}$, as specified below:

$$b_{L_{max}} = x_1 - \sqrt{x_1^2 - n}, \qquad (10.5)$$

i.e. the adopted value of $b_L$ must always satisfy the following condition:

$$b_L \leq b_{L_{max}} = x_1 - \sqrt{x_1^2 - n}, \qquad (10.6)$$

since, otherwise, the equation that defines the value of $i_L$, which will be presented below, is no longer applicable. Therefore, before proceeding to calculate the value of $i_L$, it is necessary to verify that the value of $b_L$ to be adopted satisfies the condition $b_L \leq b_{L_{max}}$.

That said, the value of $i_L$, including the iteration for $x_1$, is given by the following general equation:

$$i_L = \left\lfloor \frac{\frac{a_L + b_L}{2} - x_1}{s} \right\rfloor + 1, \qquad (10.7a)$$

i.e.:

$$i_L = \left\lfloor \frac{\frac{a_L + b_L}{2} - \frac{n - \left\lfloor \frac{n - 2 \cdot \lfloor \sqrt{n} \rfloor}{u \cdot s} \right\rfloor \cdot u \cdot s + 1}{2}}{s} \right\rfloor + 1, \qquad (10.7b)$$

which is applicable for any step-size factor $s$ (but, keeping always in mind that, by adopting a step-size factor $s$ greater than 2, we will never be certain of the primality of $n$, even though we reach the limit number of iterations $i_L$without having identified any pair of non-trivial factors of $n$), and applied here for $s = 2$:

$$i_L = \left\lfloor \frac{\frac{a_L + b_L}{2} - x_1}{s} \right\rfloor + 1 = \left\lfloor \frac{\frac{a_L + b_L}{2} - x_1}{2} \right\rfloor + 1, \qquad (10.8a)$$

i.e.:

$$i_L = \left\lfloor \frac{\frac{a_L + b_L}{2} - \frac{n - \left\lfloor \frac{n - 2 \cdot \lfloor \sqrt{n} \rfloor}{u \cdot 2} \right\rfloor \cdot u \cdot 2 + 1}{2}}{2} \right\rfloor + 1, \qquad (10.8b)$$

denoting by $a_L$ the largest possible value to consider of $a$ non-trivial of $n$ and as determined below:

$$a_L = \text{the largest odd positive integer contained in } \frac{n}{b_L}. \qquad (10.9)$$

Thus, with $s = 2$, the value of $x$ limit, which we denote by $x_L$, will be determined as follows:

$$x_L = x_1 + s \cdot \left\lfloor \frac{\frac{a_L + b_L}{2} - x_1}{s} \right\rfloor = x_1 + 2 \cdot \left\lfloor \frac{\frac{a_L + b_L}{2} - x_1}{2} \right\rfloor, \qquad (10.10a)$$

i.e.:

$$x_L = x_1 + s \cdot (i_L - 1) = x_1 + 2 \cdot (i_L - 1). \qquad (10.10b)$$

To simplify (avoiding the a priori determination of the value of $a_L$), the value of $i_L$ can also be calculated as follows (including the iteration for $x_1$), applied here for $s = 2$ but applicable, as a general equation, for any step-size factor $s$:

$$i_L = \left\lfloor \frac{\frac{\frac{n}{b_L} + b_L}{2} - x_1}{s} \right\rfloor + 1 = \left\lfloor \frac{\frac{n + b_L^2}{2 \cdot b_L} - x_1}{2} \right\rfloor + 1. \qquad (10.11)$$

Thus, with $s = 2$, the value of $x$ limit, i.e. $x_L$, will be determined as follows:

$$x_L = x_1 + s \cdot \left\lfloor \frac{\frac{n + b_L^2}{2 \cdot b_L} - x_1}{s} \right\rfloor = x_1 + 2 \cdot \left\lfloor \frac{\frac{n + b_L^2}{2 \cdot b_L} - x_1}{2} \right\rfloor, \qquad (10.12)$$

i.e., as previously shown:

$$x_L = x_1 + s \cdot (i_L - 1) = x_1 + 2 \cdot (i_L - 1). \qquad (10.10b)$$

Obviously, for large or very large values of $n$, to calculate the value of $i_L$ we will never adopt the value of $b_L = 3$ (which, in such cases, will still certainly satisfy the condition $b_L \leq b_{L_{max}}$) since, with this value of $b_L$, in such cases we would have a relatively large or relatively very large value of $i_L$ (relative to the order of magnitude of $n$). Indeed, excluding the cases of small or very small values of $n$ (for which we must always verify that the value of $b_L$ to be adopted satisfies the condition $b_L \leq b_{L_{max}}$), to calculate the value of $i_L$ we will choose to adopt, as the value of $b_L$, a prime greater than 3 (while always verifying a priori that the chosen value of $b_L$ satisfies the condition $b_L \leq b_{L_{max}}$), subjecting $n$, in advance, to the Trial Division Algorithm (TDA) [5, pp. 133-134], attempting to divide it by all primes $p$ such that $2 < p \leq p_f$, denoting by $p_f$ the largest prime considered and contained in $\sqrt{n}$ but anyway relatively small (relative to the order of magnitude of $n$). So that, in case the division attempts were not successful by dividing $n$ by all primes $p$ such that $2 < p \leq p_f$, we will adopt, as the value of $b_L$, the first prime following $p_f$, so as to stop at the pair of values of $x$ and $y$ that determines the value of $x_L - y_L = b_L$ or the smallest value of $x_L - y_L$ greater than $b_L$, highlighting, where necessary, that, starting from the iteration immediately following the one that coincides with the limit number of iterations $i_L$, in both cases we will always have the condition $x_L - y_L < b_L$, a condition which, as we already know, is excluded by having previously subjected $n$ to the Trial Division Algorithm (TDA), attempting to divide it by all primes $p$ such that $2 < p \leq p_f$ (this, of course, also holds when we adopt the value of $b_L = 3$, i.e. when $n$ has not previously been subjected to the Trial Division Algorithm).

Still regarding the preliminary application of the Trial Division Algorithm (TDA) in order to reach, in cases of primality, the limit number of iterations with a smaller number of iterations, it is clear from (10.6) that we can decide a priori to adopt a value for $i_L$ that is equal to or close to a specific chosen value that we denote by $i_{LS}$ (thus, even with $i_{LS} = 1$). This aspect is very important for large or very large values of $n$, for which we decide to perform a larger number of division attempts with TDA, but we significantly reduce—compared to the number of division attempts with TDA—the number of iterations required to reach the limit number of iterations $i_L$. To do this (therefore, for large or very large values of $n$), we must first calculate the value of $i_L$ with $b_L = 3$, then choose a value for $i_{LS}$—which is clearly less than the value of $i_L$ calculated with $b_L = 3$—and then adopt, as the value of $b_L$ to calculate the value of $i_L$, the largest prime less than the following value, denoted by $b_{L_p}$, calculated using the chosen value of $i_{LS}$ (to standardize the procedure, since if we choose $i_{LS} = 1$ and adopt the smallest prime greater than $b_{L_p}$, the equation that defines the value of $i_L$, i.e. (10.11), is no longer applicable as we already know, and clearly we proceed to calculate the value of $i_L$ only if the attempts to divide $n$ by all primes preceding the determined value of $b_L$ have been unsuccessful):

$$b_{L_p} = x_1 + s \cdot \left(i_{LS} - 1\right) - \sqrt{\left[x_1 + s \cdot \left(i_{LS} - 1\right)\right]^2 - n}, \qquad (10.13a)$$

i.e., with $s = 2$:

$$b_{L_p} = x_1 + 2 \cdot \left(i_{LS} - 1\right) - \sqrt{\left[x_1 + 2 \cdot \left(i_{LS} - 1\right)\right]^2 - n}, \qquad (10.13b)$$

and, therefore, we adopt the following value of $b_L$:

$$b_L = \text{the largest prime less than } b_{L_p}, \qquad (10.14)$$

and consequently, we preliminarily subject $n$ to the Trial Division Algorithm (TDA), attempting to divide it by all primes preceding the value of $b_L$ thus determined, end finally (if the attempts to divide $n$ by all primes preceding the determined value of $b_L$ have been unsuccessful), we calculate the value of $i_L$ for this value of $b_L$, using (10.11).

Clearly, if a preliminary application of the Trial Division Algorithm (TDA) to $n$ were to reveal small factors of $n$, we would end up calculating the value of $i_L$ for a smaller value of $n$ than the one initially considered, for which we would then have to verify the factorization.

## 11. Example of Application of DFM-1 (with $s = 2$)

In this section, we proceed with one final, detailed example of the application of DFM-1 (with $s = 2$), also reporting the simplifications described in the last part of Section 5 of this paper, which, while always leaving the number of iterations to be performed unchanged, determine an overall reduction in the computational cost.

Let $n$ = 3986359420010593 be the semiprime, which is the product of the primes $a$ = 87281521 and $b$ = 45672433 [and with $\varphi_s(n) = (a-1)\cdot(b-1) = 3986359287056640$], of which, with the adoption of FFA-1 starting from the value of $x = x_1 = \lceil\sqrt{n}\rceil = 63137623$ as the smallest integer greater than $\sqrt{n} = 63137622.85682\ldots$, the unique pair of non-trivial factors is identified at the following number of iterations $i_{C_F}$ (including the iteration for $x_1 = \lceil\sqrt{n}\rceil$) [see (4.3)]:

$$i_{C_F} = \frac{a+b}{2} - \lfloor\sqrt{n}\rfloor = \frac{87281521 + 45672433}{2} - 63137622 = 3339355,$$

that is, with $s = 1$ (see (9.1)):

$$i_{C_F} = \frac{\frac{a+b}{2} - \lceil\sqrt{n}\rceil}{s} + 1 = \frac{\frac{87281521 + 45672433}{2} - 63137623}{1} + 1 = 3339355.$$

That said, below we proceed to identify the unique pair of non-trivial factors $a$ and $b$ of $n$ by adopting DFM-1 (with $s = 2$), but we do not proceed, in advance, both to verify that the condition $x = \sqrt{n}$ does not occur and to subject $n$ to a quick probabilistic test of primality, since this is an example for which we already know that these conditions do not occur, and we choose not to resort, in advance, to the Trial Division Algorithm (TDA), thus adopting the value of $b_L = 3$ (which, in this case, clearly satisfies the condition $b_L \le b_{L_{max}}$) to determine the limit number of iterations $i_L$.

First, we calculate the initial value $x = x_1$ [see (9.2)]:

$$\begin{aligned} x_1 &= \frac{n - \varphi_{s1}(n) + 1}{2} \\ &= \frac{n - \left\lfloor \frac{\varphi_{s1_p}(n)}{u\cdot s} \right\rfloor \cdot u\cdot s + 1}{2} \\ &= \frac{n - \left\lfloor \frac{n - 2\cdot\lfloor\sqrt{n}\rfloor}{2\cdot 2} \right\rfloor \cdot 2\cdot 2 + 1}{2} \\ &= 63137623, \end{aligned}$$

with the following values of $\varphi_{s1_p}(n)$ and $\varphi_{s1}(n)$ [see (6.7) and 6.13), respectively]:

$$\varphi_{s1_p}(n) = n - 2\cdot\lfloor\sqrt{n}\rfloor = 3986359420010593 - 2\cdot 63137622 = 3986359293735349,$$

$$\varphi_{s1}(n) = \left\lfloor \frac{\varphi_{s1_p}(n)}{u\cdot s} \right\rfloor \cdot u\cdot s = \left\lfloor \frac{3986359293735349}{2\cdot 2} \right\rfloor \cdot 2\cdot 2 = 3986359293735348.$$

Therefore, we now proceed with the phases of DFM-1:

**1st phase** - Check of the values of $y = y_1$ for $x = x_1$ [see (9.3)]:

$$
\begin{aligned}
y_1 &= \sqrt{x_1^2 - n} \\
&= \sqrt{63137623^2 - 3986359420010593} \\
&= \sqrt{3986359438090129 - 3986359420010593} \\
&= \sqrt{18079536} \\
&= 4252{,}00376\ldots,
\end{aligned}
$$

and wanting to consider the simplification described in the last part of Section 5 of this paper (regarding the possible elimination of the square root operation by observing the last digit of $y_1^2$), we have:

$$y_1^2 = x_1^2 - n = 63137623^2 - 3986359420010593 = 18079536,$$

with the last digit of $y_1^2$ being 6, so, since such digit indicates that $y_1^2$ could be a perfect square, we will necessarily have to perform the square root operation of $y_1^2$.

Therefore, having unsuccessfully checked for $x = x_1$ since the corresponding value of $y = y_1$ was found to be a non-integer, i.e. since the value of $y^2 = y_1^2$ was found to be a non-perfect square, we repeat the 1st phase, increasing the value of $x$ always by two units at a time (starting, therefore, from the 2nd iteration with the value of $x_2 = x_1 + 2 = \frac{n - \left\lfloor \frac{n - 2 \cdot \lfloor \sqrt{n} \rfloor}{4} \right\rfloor \cdot 4 + 1}{2} + 2$) in the equation $y = \sqrt{x^2 - n}$ (i.e. in the equation $y^2 = x^2 - n$)—or, resorting to the simplification described in the last part of Section 5 of this paper, by adopting the equation $y_r^2 = y_f^2 + 4 \cdot (x_f + 1)$, as given in Section 5 with (5.5), for each iteration subsequent to the 1st iteration, starting, therefore, from the 2nd iteration with the value of $x_f = x_1 = \frac{n - \left\lfloor \frac{n - 2 \cdot \lfloor \sqrt{n} \rfloor}{4} \right\rfloor \cdot 4 + 1}{2}$ and, thereafter, increasing the value of $x_f$ always by two units at a time—until we have identified the integer value of $y$ with $x - y \neq 1$ (or the integer value of $y_r$ with $x_f + 2 - y_r \neq 1$, where $x_f$ is the value of $x$ corresponding to the iteration preceding the one in which we have identified the integer value of $y_r$ with $x_f + 2 - y_r \neq 1$), which is $y_C$, with which, together with the value of $x_C$, we identify, with the 2nd phase, the unique pair of non-trivial factors $a$ and $b$ of $n$ (obviously, since this is an example, we already know that only one pair of non-trivial factors of $n$ occurs). Thus, at the 1669677th iteration starting from the 2nd iteration, we identify, with the value of $x = x_C = 66476977$, the following integer value of $y$ (with $x - y \neq 1$) as $y_C$ (therefore, with $i_{C_D} = 1669678$):

$$
\begin{aligned}
y_C &= \sqrt{x_C^2 - n} \\
&= \sqrt{\left[x_1 + s \cdot \left(i_{C_D} - 1\right)\right]^2 - n} \\
&= \sqrt{[63137623 + 2 \cdot (1669678 - 1)]^2 - 3986359420010593} \\
&= \sqrt{(63137623 + 2 \cdot 1669677)^2 - 3986359420010593} \\
&= \sqrt{66476977^2 - 3986359420010593} \\
&= \sqrt{4419188471058529 - 3986359420010593} \\
&= \sqrt{432829051047936} \\
&= 20804544;
\end{aligned}
$$

or, wanting to consider the simplifications described in the last part of Section 5 of this paper, we have the value of $y_C$ calculated as indicated below, with $x_f = 66476975$ and with $y_f^2 = 432828785140032$ that, in the previous iteration, was $y_r^2$ and for which we had not performed the square root operation since, in that case, the last digit of $y_r^2$ was 2 and, therefore, certainly such $y_r^2$ could not have been a perfect square:

$$
\begin{aligned}
y_C &= \sqrt{y_r^2} \\
&= \sqrt{y_f^2 + 4 \cdot \left(x_f + 1\right)} \\
&= \sqrt{432828785140032 + 4 \cdot (66476975 + 1)} \\
&= \sqrt{432828785140032 + 4 \cdot 66476976} \\
&= \sqrt{432828785140032 + 265907904} \\
&= \sqrt{432829051047936} \\
&= 20804544,
\end{aligned}
$$

having performed the square root operation of this $y_r^2$, i.e. $y_r^2 = 432829051047936$, since the last digit of this $y_r^2$ is 6 and, therefore, it indicates that this $y_r^2$ could be a perfect square.

Therefore, with the value of $y_C$, together with the value of $x_C = 66476977$ (which, if we consider the simplification described in the last part of Section 5 of this paper by adopting (5.5), is $x_C = x_f + 2$, where $x_f = 66476975$ is the value of $x$ corresponding to the iteration preceding the one in which we have identified the integer value of $y_r = y_C$, i.e. the integer value of $y_r$ with $x_f + 2 - y_r \neq 1$, with the values of $x_f$ and $y_r$ that appear in the same equation), we identify, with the 2nd phase, the unique pair of non-trivial factors *a* and *b* of *n*:

**2nd phase** - Calculation of the values of $a$ and $b$ [see (9.6), (9.7), (9.8), and (9.9)]:

$$a = x_C + y_C = 66476977 + 20804544 = 87281521,$$

$$b = x_C - y_C = 66476977 - 20804544 = 45672433.$$

It should be noted that the value of $i_{C_D} = \frac{\frac{a+b}{2} - x_1}{s} + 1 = \frac{\frac{a+b}{2} - x_1}{2} + 1 = 1669678$ (i.e. $i_{C_D} = \frac{\varphi_{s1}(n) - \varphi_s(n)}{u \cdot s} + 1 = \frac{\varphi_{s1}(n) - \varphi_s(n)}{2 \cdot 2} + 1 = 1669678$), which is inclusive of the iteration for $x_1 = \frac{n - \left\lfloor \frac{n - 2 \cdot \lfloor \sqrt{n} \rfloor}{u \cdot s} \right\rfloor \cdot u \cdot s + 1}{2} = \frac{n - \left\lfloor \frac{n - 2 \cdot \lfloor \sqrt{n} \rfloor}{2 \cdot 2} \right\rfloor \cdot 2 \cdot 2 + 1}{2}$, is of the following form [see (7.1$b$)]:

$$i_{C_D} = \left\lfloor \frac{i_{C_F}}{s} \right\rfloor + 1 = \left\lfloor \frac{3339355}{2} \right\rfloor + 1 = \left\lceil \frac{i_{C_F}}{s} \right\rceil = \left\lceil \frac{3339355}{2} \right\rceil = 1669678.$$

At this point, admitting that we do not know a priori (since this is an example) that only one pair of non-trivial factors of $n$ occurs, we proceed to check the factorization (or possible primality) of these identified non-trivial factors $a$ and $b$ of $n$ (after checking that $a$ and $b$ are not perfect squares, also resorting to the simplification described in the last part of Section 5 of this paper regarding the elimination of the square root operation of $b$ by observing that the last digit of $b$ is 3, and after having preliminarily subjected $a$ and $b$, unsuccessfully as example, to a quick probabilistic test of primality), by adopting DFM-1 (with $s = 2$), with two new specific applications (with or without the simplifications described in the last part of Section 5 of this paper), and determining the primality of $a$ and $b$, reaching, respectively, the limit numbers of iterations $i_{L_a}$ and $i_{L_b}$ (with $b_L = 3$, which clearly satisfies the condition $b_L \le b_{L_{max}}$ for both $a$ and $b$), shown below [see (10.11)], without having identified any pair of non-trivial factors (and with the *values* of $i_{L_a}$ and $i_{L_b}$ which always include, respectively, the iteration for $x_{1_a}$ and $x_{1_b}$):

$$i_{L_a} = \left\lfloor \frac{\frac{\frac{a}{b_L} + b_L}{2} - x_{1_a}}{s} \right\rfloor + 1 = \left\lfloor \frac{\frac{a + {b_L}^2}{2 \cdot b_L} - x_{1_a}}{s} \right\rfloor + 1 = \left\lfloor \frac{\frac{a+9}{6} - 9343}{2} \right\rfloor + 1 = 7268790,$$

where:

$$x_{1_a} = \frac{a - \left\lfloor \frac{a - 2 \cdot \lfloor \sqrt{a} \rfloor}{u \cdot s} \right\rfloor \cdot u \cdot s + 1}{2} = \frac{a - \left\lfloor \frac{a - 2 \cdot \lfloor \sqrt{a} \rfloor}{2 \cdot 2} \right\rfloor \cdot 2 \cdot 2 + 1}{2} = 9343,$$

with the following values of $x_{L_a}$ [see (10.10$b$)], $y_{L_a}$, and $x_{L_a} - y_{L_a}$:

$$x_{L_a} = x_{1_a} + s \cdot \left(i_{L_a} - 1\right) = 9343 + 2 \cdot (7268790 - 1) = 14546921,$$

$$y_{L_a} = \sqrt{x_{L_a}^2 - a} = \sqrt{14546921^2 - 87281521} = 14546917.9999998625\ldots,$$

and

$$x_{L_a} - y_{L_a} = 3.0000001374\ldots$$

as the smallest value of $x - y$ greater than $b_L = 3$;

$$i_{L_b} = \left\lfloor \frac{\frac{\frac{b}{b_L} + b_L}{2} - x_{1_b}}{s} \right\rfloor + 1 = \left\lfloor \frac{\frac{b + b_L{}^2}{2 \cdot b_L} - x_{1_b}}{s} \right\rfloor + 1 = \left\lfloor \frac{\frac{b+9}{6} - 6759}{2} \right\rfloor + 1 = 3802658,$$

where:

$$x_{1_b} = \frac{b - \left\lfloor \frac{b - 2 \cdot \lfloor \sqrt{b} \rfloor}{u \cdot s} \right\rfloor \cdot u \cdot s + 1}{2} = \frac{b - \left\lfloor \frac{b - 2 \cdot \lfloor \sqrt{b} \rfloor}{2 \cdot 2} \right\rfloor \cdot 2 \cdot 2 + 1}{2} = 6759,$$

with the following values of $x_{L_b}$ [see (10.10*b*)], $y_{L_b}$, and $x_{L_b} - y_{L_b}$:

$$x_{L_b} = x_{1_b} + s \cdot (i_{L_b} - 1) = 6759 + 2 \cdot (3802658 - 1) = 7612073,$$

$$y_{L_b} = \sqrt{x_{L_b}^2 - b} = \sqrt{7612073^2 - 45672433} = 7612069.9999997372\ldots,$$

and

$$x_{L_b} - y_{L_b} = 3.0000002627\ldots$$

as the smallest value of $x - y$ greater than $b_L = 3$.

Overall, we performed 12741126 iterations to identify the unique pair of non-trivial factors $a$ and $b$ of $n$ and to determine the primality of $a$ and $b$ (with $b_L = 3$).

While, if we had decided to continue, with the first application of DFM-1 (with $s = 2$) with or without the *simplifications* described in the last part of Section 5 of this paper, in the iterations to identify any further pairs of non-trivial factors of $n$, we would have performed, overall, 332196586765406 iterations to reach the limit number of iterations $i_L$ (always with $b_L = 3$) [see (10.11)], shown below, and be certain that further pairs of non-trivial factors of $n$, which are subsequent to the one previously identified, can no longer occur (and with the value of $i_L$ which always includes the iteration for $x_1$):

$$i_L = \left\lfloor \frac{\frac{\frac{n}{b_L} + b_L}{2} - x_1}{s} \right\rfloor + 1$$

$$= \left\lfloor \frac{\frac{n + b_L{}^2}{2 \cdot b_L} - x_1}{s} \right\rfloor + 1$$

$$= \left\lfloor \frac{\frac{n+9}{6} - 63137623}{2} \right\rfloor + 1$$

$$= 332196586765406,$$

where:

$$x_1 = \frac{n - \varphi_{s1}(n) + 1}{2}$$

$$= \frac{n - \left\lfloor \frac{\varphi_{s1_p}(n)}{u \cdot s} \right\rfloor \cdot u \cdot s + 1}{2}$$

$$= \frac{n - \left\lfloor \frac{n - 2 \cdot \lfloor \sqrt{n} \rfloor}{2 \cdot 2} \right\rfloor \cdot 2 \cdot 2 + 1}{2}$$

$$= 63137623,$$

with the following values of $x_L$ [see (10.10*b*)], $y_L$, and $x_L - y_L$:

$$x_L = x_1 + s \cdot (i_L - 1)$$

$$= 63137623 + 2 \cdot (332196586765406 - 1)$$

$$= 664393236668433,$$

$$y_L = \sqrt{x_L{}^2 - n}$$

$$= \sqrt{664393236668433^2 - 3986359420010593}$$

$$= 664393236668429{,}9999999\ldots,$$

and

$$x_L - y_L = 3.000000000000003\ldots$$

as the smallest value of $x - y$ greater than $b_L = 3$.

As has already been highlighted in Sections 5 and 6 of this paper, it is clear that it is always more convenient to stop at identifying the first (which could also be the only one) integer value of $y$ with $x - y \neq 1$ (or integer value of $y_r$ with $x_f + 2 - y_r \neq 1$), provided that there exists an integer value of $y$ with $x - y \neq 1$ (or an integer value of $y_r$ with $x_f + 2 - y_r \neq 1$), proceeding to check the factorization (or possible primality), by adopting DFM-1 again (with $s = 2$) with or without the *simplifications* described in the last part of Section 5 of this paper, of the first (or unique, but we cannot know a priori) identified non-trivial factors $a$ and $b$ of $n$ (obviously, after checking that $a$ and $b$ are not perfect squares, possibly also resorting to the simplification described in the last part of Section 5 of this paper regarding the possible elimination of the square root operation of $a$ and $b$ by observing the last digit of $a$ and $b$, and after having preliminarily subjected $a$ and $b$ to a quick probabilistic test of primality).

Furthermore, we already know from Section 10 that, if we had chosen a prime greater than 3 as the value of $b_L$ (after having previously and uselessly subjected the numbers being checked for factorization or primality to the Trial Division Algorithm (TDA), attempting to divide them by all primes $p$ such that $2 < p \leq p_f$, with $p_f$ as the largest prime considered and contained in $\sqrt{n}$ but anyway relatively small considering the order of magnitude of $n$), we would have obtained a reduction, in all the cases represented above, of the limit number of iterations determined above with respect to $b_L = 3$.

Thus, always admitting that we do not know a priori (since this is an example) that only one pair of non-trivial factors of $n$ occurs, if we had proceeded to check of the factorization or possible primality—always adopting DFM-1, i.e. with $s = 2$, with two new specific applications, with or without the simplifications described in the last part of Section 5 of this paper—of the non-trivial factors $a$ and $b$ of $n$ by adopting, for example and for both $a$ and $b$ (to simplify), the value of $b_L = 2543$ as the 371st prime starting from the number 3 (after verifying, for both $a$ and $b$, that this chosen value of $b_L$ satisfies the condition $b_L \leq b_{L_{max}}$, and after having previously and uselessly subjected the non-trivial factors $a$ and $b$ to the Trial Division Algorithm (TDA), attempting to divide them by all primes $p$ such that $2 < p \leq p_f$, being $p_f = 2539$ as the 370th prime starting from the number 3), we would have determined the primality of $a$ and $b$, reaching, respectively, the limit numbers of iterations $i_{L_a}$ and $i_{L_b}$, shown below, without having identified any pair of non-trivial factors:

$$i_{L_a} = \left\lfloor \frac{\frac{\frac{a}{b_L} + b_L}{2} - x_{1_a}}{s} \right\rfloor + 1 = \left\lfloor \frac{\frac{a + b_L^2}{2 \cdot b_L} - x_{1_a}}{s} \right\rfloor + 1 = \left\lfloor \frac{\frac{a + 6466849}{5086} - 9343}{2} \right\rfloor + 1 = 4545,$$

where:

$$x_{1_a} = \frac{a - \left\lfloor \frac{a - 2 \cdot \lfloor \sqrt{a} \rfloor}{u \cdot s} \right\rfloor \cdot u \cdot s + 1}{2} = \frac{a - \left\lfloor \frac{a - 2 \cdot \lfloor \sqrt{a} \rfloor}{2 \cdot 2} \right\rfloor \cdot 2 \cdot 2 + 1}{2} = 9343,$$

with the following values of $x_{L_a}$, $y_{L_a}$, and $x_{L_a} - y_{L_a}$:

$$x_{L_a} = x_{1_a} + s \cdot \left(i_{L_a} - 1\right) = 9343 + 2 \cdot (4545 - 1) = 18431,$$

$$y_{L_a} = \sqrt{{x_{L_a}}^2 - a} = \sqrt{18431^2 - 87281521} = 15887.7386685456\ldots,$$

and

$$x_{L_a} - y_{L_a} = 2543.2613314543\ldots$$

as the smallest value of $x - y$ greater than $b_L = 2543$;

$$i_{L_b} = \left\lfloor \frac{\frac{\frac{b}{b_L} + b_L}{2} - x_{1_b}}{s} \right\rfloor + 1 = \left\lfloor \frac{\frac{b + {b_L}^2}{2 \cdot b_L} - x_{1_b}}{s} \right\rfloor + 1 = \left\lfloor \frac{\frac{b + 6466849}{5086} - 6759}{2} \right\rfloor + 1 = 1747,$$

where:

$$x_{1_b} = \frac{b - \left\lfloor \frac{b - 2 \cdot \lfloor\sqrt{b}\rfloor}{u \cdot s} \right\rfloor \cdot u \cdot s + 1}{2} = \frac{b - \left\lfloor \frac{b - 2 \cdot \lfloor\sqrt{b}\rfloor}{2 \cdot 2} \right\rfloor \cdot 2 \cdot 2 + 1}{2} = 6759,$$

with the following values of $x_{L_b}$, $y_{L_b}$, and $x_{L_b} - y_{L_b}$:

$$x_{L_b} = x_{1_b} + s \cdot \left(i_{L_b} - 1\right) = 6759 + 2 \cdot (1747 - 1) = 10251,$$

$$y_{L_b} = \sqrt{{x_{L_b}}^2 - b} = \sqrt{10251^2 - 45672433} = 7707.8251147778\ldots,$$

and

$$x_{L_b} - y_{L_b} = 2543.1748852221\ldots$$

as the smallest value of $x - y$ greater than $b_L = 2543$.

Overall (checking $a$ and $b$), with the adoption of $b_L = 2543$ and proceeding, in advance, to perform only 370 attempts to divide both $a$ and $b$ with the Trial Division Algorithm (TDA), we performed 6292 iterations to determine the primality of $a$ and $b$, i.e. 11065156 iterations less than those performed overall, always for the same determination of the primality of $a$ and $b$, with the adoption of $b_L = 3$. So and as has already been highlighted in Section 6 of this paper, it is clear that it is always more convenient (as for the adoption of FFA-1) to resort, in advance, to the Trial Division Algorithm (TDA) )—possibly also adopting (10.13*b*) and (10.14) a priori—in order to reach, in cases of primality, the limit number of iterations with a smaller number of iterations. in order to reach, in cases of primality, the limit number of iterations with a smaller number of iterations.

## 12. Conclusion

In this paper, we have presented, in a very detailed way and in addition to other aspects such as the hypothetical scenario, the new Fermat-type factorization algorithm (DFM-1) with $s = 2$, which we arrived at by applying Euler's totient function, in a particular way (by modifying, for factorization purposes, Euler's totient function $\varphi(n)$ of any $n$ that is a semiprime), to FFA-1, which is referred to in this same paper for comparison purposes. The result is a halving (possibly rounding up to the next integer) of the number of iterations required by FFA-1. And this, obviously, implies a reduction—which can even be considerable—in computation time compared to that required by FFA-1 (with the reduction which, obviously, certainly does not exist in cases in which the value of the number of iterations, required by FFA-1, is equal to 1).

Therefore, we conclude this paper with the hope that the presented new factorization algorithm (DFM-1), also considering what has been described regarding the hypothetical scenario, can be combined with other integer factorization algorithms in order to achieve more significant results.